\documentclass{amsart}
\usepackage{latexsym}
\title{C$^{*}$-Structure and K-theory of Boutet de Monvel's algebra}
\author{S.T.Melo \and R.Nest \and E.Schrohe}
\date{}
\newtheorem{thm}{Theorem}
\newtheorem{pro}{Proposition}
\newtheorem{lem}[pro]{Lemma}
\newtheorem{cor}[pro]{Corollary}

\begin{document}
%
%
%
%
\newcommand{\op}{operator}
\newcommand{\ops}{operators}
\newcommand{\psd}{pseudo\-dif\-fer\-en\-tial}
\newcommand{\Psd}{Pseudo\-dif\-fer\-en\-tial}
\newcommand{\tp}{transmission property}
\newcommand{\bdm}{Boutet de Monvel}
\newcommand{\pf}{{\em Proof}: }
\newcommand{\ifoi}{if and only if}
\newcommand{\cst}{C$^{*}$}
%
%
%
%
\newcommand{\z}{{\mathbb Z}}
\newcommand{\co}{{\mathbb C}}
\newcommand{\re}{{\mathbb R}}
\newcommand{\so}{{\mathbb S}^1}
\newcommand{\skp}{{\mathbb S}^{k+1}}
\newcommand{\tk}{{\mathbb T}^k}
\newcommand{\rp}{\overline{{\mathbb R}}_{+}}
\newcommand{\rn}{{\mathbb R}^{n}}
\newcommand{\rnm}{{\mathbb R}^{n-1}}
\newcommand{\rnp}{\bar{\mathbb R}^{n}_{+}}
\newcommand{\rtd}{{\mathbb R}^{2d}}
\newcommand{\sx}{S^{*}X}
\newcommand{\six}{S^{*}\dot X}
\newcommand{\sbx}{S^{*}\partial X}
\newcommand{\sxt}{\sx/\!\!\sim}
\newcommand{\bx}{\partial X}
\newcommand{\ix}{\dot X}
\newcommand{\bw}{\partial W^k}
\newcommand{\iw}{\dot W^k}
\newcommand{\siw}{S^{*}\dot W^k}
\newcommand{\tix}{T^*\dot X}
\newcommand{\bix}{B^*\dot X}
%
%
%
%
\newcommand{\cix}{C^{\infty}(X)}
\newcommand{\ciix}{C^{\infty}(\ix)}
\newcommand{\cciix}{C_{c}^{\infty}(\ix)}
\newcommand{\cisx}{C^{\infty}(\sx)}
\newcommand{\cibx}{C^{\infty}(\bx)}
\newcommand{\cisbx}{C^{\infty}(\sbx)}
\newcommand{\cosx}{C_{0}(\six)}
\newcommand{\cotx}{C_{0}(\tix)}
\newcommand{\cosw}{C_{0}(\siw)}
\newcommand{\cobx}{C_{0}(\bix)}
\newcommand{\txix}{T^*\dot X}
\newcommand{\h}{{\mathfrak H}}
\newcommand{\lh}{{\mathfrak L}(\h)}
\newcommand{\as}{{\mathcal A}}
\newcommand{\gs}{{\mathcal G}}
\newcommand{\is}{{\mathcal I}}
\newcommand{\xs}{{\mathcal X}}
\newcommand{\ac}{{\mathfrak A}}
\newcommand{\gc}{{\mathfrak G}}
\newcommand{\kc}{{\mathfrak K}}
\newcommand{\wc}{{\mathfrak W}}
\newcommand{\ic}{{\mathfrak I}}
\newcommand{\tc}{{\mathfrak T}}
\newcommand{\xc}{{\mathfrak X}}
\newcommand{\rc}{{\mathfrak L}(L^{2}(\rp)\oplus\co)}
\newcommand{\srp}{{\mathcal S}(\rp)}
\newcommand{\ef}{_{_{EF}}}
%
%
%
%
\newcommand{\sig}{\bar\sigma}
\newcommand{\gam}{\bar\gamma}
\newcommand{\img}{\text{Im}\,\gam}
\newcommand{\ind}{\text{{\tt ind}}\,}
%
%
%
%
\newcommand{\cvb}{x^\prime,\xi^\prime}
\newcommand{\cqd}{\hfill$\Box$}
\newcommand{\supp}{{\rm supp}}
\newcommand{\pp}{P_{_{+}}}
\newcommand{\qp}{Q_{_{+}}}
\newcommand{\du}{d\!u}
\newcommand{\dv}{d\!v}
\newcommand{\dx}{d\!x}
\newcommand{\dy}{d\!y}
\newcommand{\dS}{d\!S}
\newcommand{\dxi}{d\!\xi}

\begin{abstract}

We consider the norm closure $\ac$\ of the algebra of all \ops\ of order and class zero in \bdm's calculus on a manifold $X$\ with boundary $\bx$. 
We first describe the image and the kernel of the continuous extension of the boundary principal symbol homomorphism to $\ac$. 
If $X$\ is connected and $\bx$\ is not empty, we then show that the K-groups of $\ac$\ are topologically determined. In case $\bx$\ has torsion free K-theory, we get $K_i(\ac/\kc)\simeq K_i(C(X))\oplus K_{1-i}(\cotx)$, $i=0,1$,
with $\kc$\ denoting the compact ideal, and $\tix$\ denoting the cotangent bundle of the interior.
Using \bdm's index theorem, we also prove that the above formula holds for $i=1$\ even without this torsion-free hypothesis; and show, moreover, that
$K_1(\ac)\simeq K_1(C(X))\oplus\ker\chi$, with $\chi:K_0(\tix)\to\z$\ denoting the topological index. 
For the case of orientable, two-dimensional $X$, 
$K_0(\ac)\simeq\z^{2g+m}$\ and $K_1(\ac)\simeq\z^{2g+m-1}$, where $g$\ is the genus of $X$\ and $m$\ is the 
number of connected components of $\bx$. We also obtain a composition sequence $0\subset\kc\subset\gc\subset\ac$, with $\ac/\gc$\ commutative and $\gc/\kc$\ 
isomorphic to the algebra of all continuous functions on the cosphere bundle of $\bx$\ with values in  compact operators on $L^2(\rp)$.

\end{abstract}
\maketitle

\section{Introduction}\label{intr}

\subsection{\bdm's algebra}
Let $X$\ denote an $n$-dimensional compact manifold, with boundary $\bx$, embedded in a closed manifold $\Omega$\ of the same dimension. Given a \psd\ \op\ 
$P$\ on $\Omega$\ and $u$\
in $\cix$, we denote by $\pp u$\ the restriction to the interior of $X$, $\ix$, of $P$ applied to the extension by zero of $u$\ to $\Omega$. This gives a 
continuous mapping $\pp:\cix\to\ciix$, completely determined by the restriction of $P$\ to $\ix$. It is said that $P$\ has the {\em \tp}\ 
if the image of $\pp$\ is contained in $\cix$. 

\bdm\ characterized \cite{B1}\ the \tp\ for {\em classical} (i.e., with polyhomogeneous symbols)
\psd\ \ops\ in terms of certain symmetry conditions satisfied at the boundary by the homogeneous components of the symbol and their derivatives. 
In \cite{B2}, he constructed an algebra containing all classical boundary-value problems on $X$, together with
their parametrices, in the elliptic case. The elements of his calculus, called {\em Green operators}, are matrices of the form
\begin{equation}
\label{boutet}
A=\left(\begin{array}{cc}\pp+G&K\\T&S\end{array}\right):\begin{array}{c}\cix\\\oplus\\\cibx\end{array}\longrightarrow\begin{array}{c}\cix\\\oplus\\\cibx
\end{array},
\end{equation}
where $P$\ is a \psd\ \op\ with the \tp\ on $X$\ and $S$\ is a \psd\ \op\ on the closed manifold $\bx$. The \ops\ $G$, $K$, and $T$\
are regularizing in the interior of $X$\ and locally at the boundary are given as \psd\ \ops\ on $\rnm$\ with symbols taking values in operators from
$\srp$\ to $\srp$, from $\co$\ to $\srp$, and from $\srp$\ to $\co$, respectively, where $\srp$\ denotes the restriction to the non-negative half-axis
$\rp$\ of functions in the Schwartz class ${\mathcal S}(\re)$. They are called, respectively, singular Green \ops, Poisson \ops\ and trace \ops. 

Being \psd\ \ops\ with \op-valued symbols, Green \ops\ have an order assigned to them. 
Moreover, they also have a {\em class}\ (or {\em type}), related to the order of the derivatives appearing in the 
boundary condition. There exist order and class reducing \ops, which are isomorphisms between suitable Sobolev spaces. For many purposes it is therefore 
enough to consider the \ops\ of order and class zero. Detailed expositions of \bdm's calculus can be found in \cite{G,RS,S3}. The precise estimates satisfied 
by the symbols defining the \ops\ $G$, $K$\ and $T$\ are listed in \cite{G}, Definition 2.3.13. One peculiar aspect worth to be mentioned is that, in the 
polyhomogeneous case, and for \ops\ of order $d$, the degree of homogeneity of the leading term in the asymptotic expansion for the symbols of $G$\ and $K$\ 
is $(d-1)$; while that for $P$, $T$\ and $S$\ is, as expected, $d$.

The set $\as$\ of all classical, or polyhomogeneous, Green \ops\ of order zero and class zero on $X$\ is an 
adjoint-invariant sub-algebra of $\lh$, the bounded \ops\ on the Hilbert space $\h=L^2(X)\oplus H^{-\frac{1}{2}}(\bx)$, $H^{-\frac{1}{2}}$\ denoting the 
usual Sobolev space. It is, in fact, a $\Psi^{*}$-algebra, in the sense of Gramsch \cite{Gr}\ (by \cite{S2}, Corollary~4.11; the spectral invariance of
$\as$\ had been shown earlier by Schulze \cite{Sch}). It then follows that $\as$\ is invariant under the holomorphic functional calculus and its K-theory 
coincides with that of its norm closure, which we denote by $\ac$.

There is some degree of arbitrariness with respect to how one defines the order of Poisson and trace \ops. We adopt Grubb's definitions in order to be able 
to freely quote from her book. Besides, that coincides with \bdm's original definitions. On the other hand, since we are interested only in $L^2$-theory,
it would also be natural to say that a trace or a Poisson \op\ is 
of order zero if its \op-valued symbol satisfies zero-order norm estimates. If we did that, not only we would get rid of some uncomfortable $\pm 1/2$'s in 
the exponent (\cite{G}, (1.2.19) and (2.3.47), for example), but also the zero-order Green \ops\ would form an algebra of bounded \ops\ on the more familiar
Hilbert space $L^2(X)\oplus L^2(\bx)$. The two approaches are, of course, equivalent. One could go from one algebra to the other by
\[
\as\ni\left(\begin{array}{cc}\pp+G&K\\T&S\end{array}\right)
\longmapsto
\left(\begin{array}{cc}1&0\\0&\Lambda^{\frac{1}{2}}\end{array}\right)
\left(\begin{array}{cc}\pp+G&K\\T&S\end{array}\right)
\left(\begin{array}{cc}1&0\\0&\Lambda^{-\frac{1}{2}}\end{array}\right),
\]
where $\Lambda=(1-\triangle_{\bx})^{-\frac{1}{2}}$, $\triangle_{\bx}$\ denoting a second-order nonpositive elliptic operator on $\bx$.

\subsection{The boundary principal symbol}
Two homomorphisms are defined on the algebra of all classical Green \ops\ on $X$ (\cite{B2}, \S 4): the principal 
symbol, which we denote by $\sigma$; and the boundary principal symbol, which we denote by $\gamma$. Given $A$\ as in (\ref{boutet}), $\sigma(A)=\sigma(P)$\ 
is simply the usual principal symbol of $P$, regarded as a function on the cosphere bundle $\sx$, which we consider for convenience as the set of all unit 
covectors with respect to a chosen Riemannian metric on $X$. Since the singular Green operators are regularizing in the interior of $X$, 
$\sigma(P)=\sigma(P^\prime)$\ when $\pp+G=\pp^\prime+G^\prime$. 

The boundary principal symbol of $A$, $\gamma(A),$\ depends smoothly 
on covectors in $\sbx$, takes values in Green \ops\ on $\rp,$\ and needs local coordinates to be described. Let $p(x^\prime,x_n,\xi^\prime,\xi_n)$, 
$g(x^\prime,\xi^\prime,\xi_n,\eta_n)$, $k(x^\prime,\xi^\prime,\xi_n)$, $t(x^\prime,\xi^\prime,\xi_n),$\ and $s(x^\prime,\xi^\prime)$\ be the symbols
of $P$, $G$, $K$, $T$, and $S$ with respect to coordinates $x=(x^\prime,x_n)$\ for which the boundary is 
\label{coord}
$\{x_n=0\}$. Let $p_0$, $g_0$, $k_0$, $t_0$, and
$s_0$\ denote the leading terms in the asymptotics expansions of $p$, $g$, $k$, $t$\ and $s$. For each $(x^\prime,\xi^\prime)$, one then defines 
$\gamma(A)(x^\prime,\xi^\prime)$\ as the Green \op\ on $\rp$\ obtained from the symbols (regarded as functions of $\xi_n$\ and $\eta_n$) 
$p_0(x^\prime,0,\xi^\prime,\xi_n)$, $g_0(x^\prime,\xi^\prime,\xi_n,\eta_n)$, $k_0(x^\prime,\xi^\prime,\xi_n)$, $t_0(x^\prime,\xi^\prime,\xi_n)$, and 
$s_0(x^\prime,\xi^\prime)$. 
\begin{equation}
\label{pbs}
\gamma(A)(\cvb)=\left(\begin{array}{cc}p_0(x^\prime,0,\xi^\prime,D_n)_{_{+}}+g_0(\cvb,D_n)&k_0(\cvb,D_n)\\t_0(\cvb,D_n)&s_0(\cvb)\end{array}\right)
\end{equation}
A word about the notation: The singular Green operator 
$G$ acts like a pseudodifferential operator along
the boundary, taking values in regularizing operators in the normal direction.
It has a representation as an operator with a so-called symbol-kernel
$\tilde g = \tilde g(x',\xi',x_n,y_n)$, which is a function in 
${\mathcal S}(\rp \times \rp)$  for fixed $(x',\xi')$ and satisfies 
special estimates, combining the usual pseudodifferential
estimates in $x'$ and $\xi'$ with those for rapidly decreasing functions in
$x_n$ and $y_n$. 
The symbol $g$ is defined by Fourier and inverse Fourier transform:
\begin{equation}\label{symbolkernel}
g(x',\xi',\xi_n,\eta_n) = F_{x_n\to \xi_n}\overline F_{y_n\to \eta_n}\tilde
g(x',\xi',x_n,y_n).
\end{equation}
It has an expansion into homogeneous terms; the leading one is $g_0$. Inverting
the operation above, we can associate with $g_0$ a symbol-kernel 
$\tilde g_0(x',\xi',x_n,y_n)$ which is rapidly decreasing in $x_n$ and $y_n$
for fixed $(x',\xi')$. 
One denotes by $g_0(x',\xi',D_n)$ the (compact) operator
induced on $L^2(\rp)$ by this
kernel. 
Similarly,  $K$ and $T$ have symbol-kernels $\tilde k(x',\xi',x_n)$ and 
$\tilde t(x',\xi',y_n)$; these are  rapidly decreasing 
functions for fixed $(\cvb)$. 
The symbols $k$ and $t$ are defined as their Fourier and inverse 
Fourier transforms. They have asymptotic expansions with leading terms
$k_0$ and $t_0$. Via the symbol-kernels $\tilde k_0$ and $\tilde t_0$
one defines $k_0(\cvb, D_n):\co \to L^2(\rp)$ as multiplication 
by $\tilde k_0(\cvb,\cdot)$, while 
$t_0(\cvb,D_n):L^2(\rp)\to \co$ is the
operator $\varphi \mapsto \int \tilde t_0(\cvb,\cdot)\varphi$.   
Details can be found in \cite {G}, Sections 1.2 and 2.3. 

For the invariance of the above definition, see \cite{RS}, 2.3.3.1, Theorem~3; and \cite{G}, Theorem~2.4.11. Actually, the homomorphism $\gamma$\ depends
on the choice of a normal coordinate, even though the algebra and the principal symbol do not. So, let us assume that such a choice has been made; or 
equivalently, let us suppose that all changes among the above described coordinates preserve $x_n$.

On $\as$, $\sigma$\ and $\gamma$\ are *-homomorphisms. Moreover, since the Green \ops\ of order and class zero on $\rp$\ are bounded, one has 
$\gamma(A)\in C^{\infty}(\sbx,\rc)$, for all $A\in\as$.

Gohberg \cite{Goh}\ and Seeley \cite{See}\ established the equality between the norm, modulo compacts, of a singular integral \op\ on a compact manifold and 
the supremum norm of its symbol. Proofs of that estimate for \psd\ \ops\ appeared in \cite{Hp,KN}.  
We need the following generalization:
\begin{equation}
\label{rs}
\inf_{C\in\kc}||A+C||=\max\{||\sigma(A)||,||\gamma(A)||\}, \ \text{for all}\ A\in\as,
\end{equation}
with $\kc$\ denoting the ideal of the compact \ops\ on $\h$, $||\sigma(A)||$ the supremum norm of $\sigma(A)$\ on $\sx$, and $||\gamma(A)||$ the supremum
over all $(\cvb)$\ in $\sbx$\ of $||\gamma(A)(\cvb)||_{\rc}$. 
This result can be found in Rempel and Schulze's book (\cite{RS}, 2.3.4.4, Theorem~1); they credit Grubb and Geymonat \cite{GG}\ for earlier work.
Theorem~3.23 in \cite{rsm}\ further generalizes (\ref{rs}) to a larger algebra, on which a third symbol is defined. 

One obtains from (\ref{rs}), in particular, that $\sigma$\ and $\gamma$\ can be
extended to \cst-algebra homomorphisms, $\sig:\ac\to C(\sx)$\ and $\gam:\ac\to C(\sbx)\otimes\rc$, with
\begin{equation}
\label{rsc}
\inf_{C\in\kc}||A+C||=\max\{||\sig(A)||,||\gam(A)||\}, \ \text{for all}\ A\in\ac .
\end{equation}
We have written $\otimes$\ for the \cst-algebra tensor product, noting that $C(\sbx)$\ is nuclear. Equivalently, for each $A\in\ac$, $\gam(A)$\ is 
a continuous function on the cosphere bundle of the boundary, with values in bounded \ops\ on $L^{2}(\rp)\oplus\co$. 

\subsection{Statement of results and definitions}
Let $\is$\ denote the algebra of all Green \ops 
\begin{equation}
\label{dis}
\left(\begin{array}{cc}\varphi P\psi+G&K\\T&S\end{array}\right),
\end{equation} 
where $P$\ is a zero-order classical \psd\ \op\ on $X$; $G$, $K$, $T$,
and $S$\ have negative order; and $\varphi$\ and $\psi$\ belong to $\cciix$, the space 
of smooth functions with support contained in $\ix$\ (we denote by the same symbols also the operators of multiplication by $\varphi$\ or $\psi$). In 
Section~\ref{kg}, we prove that the kernel of $\gam$\ is equal to the norm closure of $\is$, which we denote by $\ic$. The crucial step for that is a 
norm estimate, modulo $\ic$, stated in Lemma~\ref{crucial}; which is, in its essence, a result for manifolds {\em with}\ boundary, 
in the sense that it gives trivial information when applied for a manifold without boundary (for then $\ic=\ac$). The usefulness of this description of 
$\ker\gam$\ to our K-theoretic calculations follows from the fact that the quotient $\ic/\kc$\ is isomorphic to the algebra $\cosx$\ of all continuous 
functions on the cosphere bundle of the interior, $\six$, that vanish at the boundary. 

This characterization of $\ker\gam$\ is equivalent to saying that, if one enlarges the ideal on the left-hand side of (\ref{rs}), then the principal symbol 
is not needed on the right. More precisely, we have (Corollary~\ref{ami}):
\[
\inf_{A^\prime\in\ic}||A+A^\prime||\,=\,||\gamma(A)||, \ \text{for all}\ A\in\as .
\]

Let $\so=\{z\in\co;\,|z|=1\}$\ and let $U:L^2(\so)\to L^2(\re)$\ denote the unitary mapping
\begin{equation}
\label{unitary}
Ug(t)=\frac{\sqrt{2}}{1+it}\, g\!\left(\frac{1-it}{1+it}\right).
\end{equation}
By ${\mathcal H}_{-1}$\ denote the image of $C^{\infty}(\so)$\ under $U$, and by ${\mathcal H}_{0}={\mathcal H}_{-1}\oplus\co$\ the direct sum of
${\mathcal H}_{-1}$\ with the constant functions. If $p(x^\prime,x_n,\xi^\prime,\xi_n)$\ is the symbol of a classical zero-order \psd\ \op\ with 
the \tp\ on $X$, with respect to local coordinates as above, then, for each $(\cvb)$, $p(x^\prime,0,\xi^\prime,\cdot)$\ belongs to ${\mathcal H}_{0}$. 
Moreover, the image by $U$\ of the Hardy space $H^2$\ is equal to $F(L^2(\rp))$, with $F$\ denoting the Fourier transform on $\re$; and the bounded \op\ 
$p(x^\prime,0,\xi^\prime,D_n)_{_{+}}$\ equals \label{pg}$F^{-1}UT_{p_{\cvb}}U^{-1}F,$\ with $T_{p_{\cvb}}$\ denoting the Toeplitz \op\ of symbol 
$p_{\cvb}(z)=p(x^\prime,0,\xi^\prime,\frac{iz-i}{z+1})$\ (we refer to \cite{Dou}\ for the definitions of Hardy space and of Toeplitz \ops). These are 
(rephrased) fundamental results for \bdm's calculus, their proofs can be found in \cite{B2}, \S 1; \cite{RS}, 2.1; and \cite{G}, 2.2. Our definition here 
of ${\mathcal H}_{-1}$\ and ${\mathcal H}_{0}$\ agrees with that of \cite{G}, but not with those of \cite{B2,RS}, where they are denoted by
${\mathcal H}_{0}$\ and ${\mathcal H}_{1}$, respectively. 

Let $\tc$\ denote the \cst-algebra of bounded \ops\ on $L^2(\rp)$\ generated by $\{p(D)_{_{+}};\,p\in{\mathcal H}_{0}\}$. The above observations prove that 
$\tc$\ is unitarily equivalent to the \cst-algebra generated by all Toeplitz \ops\ of continuous symbol; in particular, $\tc$\ contains all compact \ops\ on 
$L^2(\rp)$\ (\cite{Dou}, Proposition 7.12). 
We noted before that tbe operator $g_0(\cvb,D_n)$\ is compact. 
Hence the upper left corner of the matrix in (\ref{pbs}) belongs to $\tc$\ for every $(\cvb)$, not depending on the way we
write $\pp+G$\ as the sum of a truncated \psd\ \op\ with the \tp\ and a singular Green \op\ (see the proof of Lemma~\ref{l1}). 


Any $A\in\ac$ can be written as
$
A=\left(\begin{array}{cc}A_{11}&A_{12}\\A_{21}&A_{22}\end{array}\right),
$\ 
with the $A_{jk}$'s denoting the obvious compositions of $A$\ with orthogonal projections of $\h$\ onto its subspaces $L^2(X)$\ and $H^{-\frac{1}{2}}(\bx)$. 
By definition, $\as$\ is invariant under multiplication by these projections. We also denote $\ac_{jk}=\{A_{jk};A\in\ac\}$. Completing the matrix with zeros,
it is clear that $\ac_{11}$\ and $\ac_{22}$\ may be regarded as subalgebras of $\ac$, with non-unital inclusions; while $\ac_{12}$\ and $\ac_{21}$\ are only
subspaces. We will use the same notation for any subalgebra of $\ac$\ and a similar notation for homomorphisms. 

The preceding discussion shows that $\gam$\ maps $\ac$\ to 
\begin{equation}
\label{gama}
C(\sbx)\otimes\left(\begin{array}{cc}\tc&L^2(\rp)\\L^2(\rp)^{*}&{\mathbb C}\end{array}\right). 
\end{equation} 
In Section~\ref{ig}, we prove that the image of $\gam$, $\img$, is yet a proper subalgebra of the algebra in (\ref{gama}). The non-surjectiveness is 
observed only at the upper-left corner, as explained in the next two paragraphs. First, let us notice that every $p\in{\mathcal H}_{-1}$\ vanishes at 
infinity (that follows immediately from the definitions of $U$\ and of ${\mathcal H}_{-1}$). Hence, $p(\infty)$\ is defined for every $p$\ in 
${\mathcal H}_{0}$. 

Let $\kc$\ also denote the ideal of compact \ops\ on $L^2(\rp)$\ (we will denote by the same letter the compact 
ideal on any of the Hilbert spaces of this paper; except in the few cases when a distinction between some of them will be needed). 
It is a result of Coburn \cite{Cob1,Cob2}\ and Douglas \cite{Dou2}\ that $\tc/\kc$\ is isomorphic to 
$C(\re\cup\{\infty\})\simeq C(\so)$, with isomorphism induced by the symbol mapping (see \cite{Dou}; Theorem 7.23, and the bibliographical notes of Chapter~7,
where the influence of previous work of Gohberg is acknowledged). The mapping $p(D)_{_{+}}\mapsto p(\infty)$\ extends therefore to a $^{*}$-homomorphism 
$\lambda:\tc\to\co$\ containing $\kc$\ in its kernel, which we denote $\tc_0$. 

We show that $\img_{11}$\ contains $C(\sbx)\otimes\tc_0$, contains also $C(\bx)\otimes\co$, and that's all. The Banach-space direct sum of 
these two algebras gives $\img_{11}$; while $\gam_{jk}$\ is surjective if $j\neq 1$\ or $k\neq 1$. This description of $\img$\ is stated in 
Theorem~\ref{imagem}, in a form more suitable for applications to K-theory.
It is precisely the fact that all $\tc_0$-valued functions, but not all the $\tc$-valued ones, are contained in $\img_{11}$\ that allows our very explicit 
computation of K-groups: since $K_*(\tc_0)=0$, it follows from Theorem~\ref{imagem} that $\img$\ and $C(\bx)$\ have isomorphic K-theory 
(Corollary~\ref{kthbx}).

In Section~\ref{kt}, using that there exists a nonvanishing section of the cotangent bundle (Proposition~\ref{tritan}; for that, it is required that $X$\ is 
connected and that $\bx$\ is nonempty), we reduce to a purely topological problem the analysis of the index and exponential mappings in the six-term cyclic 
exact sequence associated to  
\begin{equation}
\label{corner}
0\to\ic/\kc\to\ac/\kc\to\ac/\ic\to 0.
\end{equation}
Theorem~\ref{kthr}\ then solves the problem of computing the K-theory of $\ac/\kc$, in the sense that both $K_0(\ac/\kc)$\
and $K_1(\ac/\kc)$\ are put in the middle of short exact sequences of abelian groups determined by the topologies of the manifold and of the cotangent
bundle of the interior, $\tix$. In case $\bx$\ has torsion-free K-theory, we get (Corollary~\ref{free}): 
\begin{equation}
\label{answer}
K_i(\ac/\kc)\simeq K_i(C(X))\oplus K_{1-i}(\cotx),\ i=0,1,
\end{equation} 
where $\cotx$\ denotes the set of all functions on $\tix$\ which get arbitrarily small outside compacts.

In Propositions~\ref{kac}\ and~\ref{kas}, we prove that the three groups $K_0(\ac/\kc)$, $K_0(\ac)$, and $K_0(\as)$\ are isomorphic; while 
$K_1(\ac)$\ and $K_1(\as)$\ are isomorphic to the kernel of the Fredholm index mapping $K_1(\ac/\kc)\to\z$, which is surjective. 

Proposition~\ref{kths}\ gives an isomorphism between $K_i(\cosx)$\ and the direct sum of $K_i(C_0(\ix))$\ and $K_{1-i}(\cotx)$, $i=0,1$. Moreover, with 
respect to this isomorphism, the canonical projection of $K_i(\cosx)$\ onto $K_{1-i}(\cotx)$\ is equivalent to the index mapping for the exact sequence 
$0\to\cotx\to\cobx\to\cosx\to 0$, with $\bix$\ denoting the bundle of unit balls over $\ix$, which may also be regarded as the radial compactfication of 
$T^*X$. 

In Section~\ref{index}, we prove (\ref{answer}) for $i=1$\ without the torsion-free assumptions of Corollary~\ref{free}. The reason why we get for $K_1$\
a better answer than for $K_0$\ is the fact, proven by \bdm~\cite{B2}, that the classical {\em difference bundle}, defined by the principal 
symbol, induces a homomorphism from $K_1(\ac/\kc)$\ to $K_0(\cotx)$, which we will denote by $\ind$. Moreover, since the difference bundle construction is 
precisely the index mapping for topological K-theory (Proposition~\ref{ruy}), the composition of $\ind$\ with the homomorphism $i_*$\ induced by the inclusion 
of $\ic/\kc$\ into $\ac/\kc$\ is equivalent to the canonical projection $K_1(C_0(\ix))\oplus K_0(\cotx)\to K_0(\cotx)$. Even though $i_*$\ is not necessarily 
injective (for the example considered in Section~\ref{superficie}, $\ker i_*\simeq\z^{m-1}$, $m\geq 1$), its restriction to $K_0(\cotx)$\ is. That is how 
$K_0(\cotx)$\ injects into $K_1(\ac/\kc)$. The injection of $K_1(C(X))$\ is induced by the embedding of $C(X)$\ as multiplication \ops. All that 
is summarized in Theorem~\ref{nosso}.

With respect to the isomorphism $K_1(\ac/\kc)\simeq K_1(C(X))\oplus K_0(\cotx)$\ of Theorem~\ref{nosso}, we show that $\ind$\ corresponds to the canonical 
projection onto $K_0(\cotx)$ (Corollary~\ref{it}). \bdm's index theorem (\cite{B2}, Theorem 5.22) then implies that $K_1(\ac)$\ is isomorphic to 
$K_1(C(X))\oplus\ker\chi$, where $\chi:K_0(\cotx)\to\z$\ denotes the topological index.

In Section~\ref{superficie}, we assume that $X$\ is a connected two-dimensional orientable manifold with nonempty boundary. We denote its genus by $g$\ and
by $m$\ the number of connected components of the boundary. We then apply the results of Section~\ref{kt}\ to prove that $K_0(\ac/\kc)$\ and $K_1(\ac/\kc)$\ 
are both isomorphic to $\z^{2g+m}$. The standard description of a closed surface as a polygon with sides identified is used to compute the K-theory of $X$\
and of $\ix$. 
 
Under our initial assumption that $X$\ is an arbitrary compact manifold with boundary, in Section~\ref{compseq}\ we give the following 
{\em composition sequence}\ 
(in the sense of \cite{D}, 4.3.2) for $\ac$: 
\begin{equation}
\label{cmpsq}
0\subset\kc\subset\gc\subset\ac, 
\end{equation}
where $\gc$\ denotes the closure of $\gs$, the algebra of all Green \ops\ $A$\ as in (\ref{boutet}) with $P$\ of negative order. All commutators of $\as$\ 
belong to $\gs$, as follows from the rules of \bdm's calculus (\cite{G}, Section 2.6, for example). In Theorem~\ref{stil}, we prove that the principal symbol 
induces an isomorphism between $\ac/\gc$\ and $C(\sxt)$, where $\sxt$\ denotes the quotient of $\sx$\ by the equivalence relation that identifies the {\em 
North and South poles}\ over each point in $\bx$\ (i.e., the two covectors that vanish on vectors tangent to the boundary). Moreover, it follows immediately 
from (\ref{rs}) that the boundary principal symbol induces an isometry on $\gc/\kc$. Then it is not hard to see (Theorem~\ref{ontogreen}) that the image of 
that isometry is equal to $C(\sbx)\otimes\kc_{\rp}$, with $\kc_{\rp}$\ denoting the ideal of compact \ops\ on $L^2(\rp)$. 

A natural problem posed by this composition sequence is to understand the connecting mappings $K_i(C({\sxt}))\to K_{1-i}(C(\sbx))$, $i=0,1$\ in the six-term 
exact sequence associated to $0\to\gc/\kc\to\ac/\kc\to\ac/\gc\to 0$. That could lead to another way of computing the K-groups of $\ac$. Moreover, as already 
suggested by the case of orientable surfaces, there may exist interesting connections between those mappings and the topology of the manifold.

\subsection{Related results}
There are many other similar composition sequences for \cst-algebras generated by \psd\ \ops. Closer to us, Cordes \cite{Cordes}\ used \cst-algebra 
techniques to solve boundary value problems on the half space. His algebra $\ac_C$\ has a composition sequence $0\subset\kc\subset{\mathfrak E}\subset\ac_C$; 
with the cosphere bundle contained in the Gelfand space of the commutative \cst-algebra $\ac_C/{\mathfrak E}$, and ${\mathfrak E}/\kc$\ isomorphic to the 
algebra of $\rp$-compact-operator-valued continuous functions on a closed subset of a compactification of the cotangent bundle of the boundary. 

For the \cst-algebra $\ac_M$\ generated by Melrose's b-\psd\ \ops\ on a compact manifold with boundary, Lauter \cite{Lauter}\ found the compostion 
sequence $0\subset\kc\subset{\mathfrak E}\subset\ac_M$, with 
$\ac_M/{\mathfrak E}$\ isomorphic (via the principal symbol) to the algebra of continuous functions on the cosphere bundle, and ${\mathfrak E}/\kc$\ 
isomorphic to a direct sum (indexed by the connected components of the boundary) of algebras of continuous $\bx$-compact-operator-valued functions on $\re$. 

For the case of a manifold with corners of dimension $n$, Melrose and Nistor \cite{MN}\ obtained a composition sequence 
$0\subset\kc=\ic_n\subset\cdots\subset\ic_1\subset\ic_0\subset\ac_M$, with $\ac_M/\ic_0$\ isomorphic to the continuous functions on the cosphere bundle, and 
each $\ic_l/\ic_{l+1}$\ isomorphic to direct sums of algebras of compact-operator-valued functions on $\re$. They explicitly computed the connecting mappings 
for the six-term exact sequences associated to each of the quotients $\ic_l/\ic_{l+1}$. Their results were generalized for algebras of \psd\ \ops\ on 
groupoids
by Monthubert \cite{Mth1,Mth2}. Also in the language of groupoids, the K-theory of $\ac_M$\ was computed by Lauter, Monthubert and Nistor \cite{LMN}, for the 
case of a manifold with connected boundary (no boundary faces of codimension larger than one).

Index theorems for \ops\ in \bdm's calculus have been established by several authors: \bdm\ \cite{B2}, Fedosov \cite{Fe}, 
Rempel and Schulze \cite{RS}, Grubb \cite{G}. The corresponding formulas, however, are not very explicit. The present project is part of our intention
to reconsider the problem under the perspective of noncommutative geometry. Computing the K-theory as well as the Hochschild and cyclic cohomology \cite{NS},
we hope to be able to adapt the index-theoretic methods developed by Nest and Tsygan \cite{NT}.

\bdm's algebra has also been studied under other \op-algebraic aspects. Fedosov, Golse, Leichtnam and Schrohe \cite{FGLS}\ showed that 
there exists a unique continuous trace on it which extends Wodzicki's noncommutative residue \cite{Wod}. Its relation to Dixmier's trace has been studied
in \cite{NSt}. Finally, it was proven by Grubb and Schrohe \cite{GS}\ that this trace is indeed given as a residue in the meromorphic extension of a suitable
\op\ trace. 

\subsection{Vector Bundles}\label{vb}
It is straightforward to generalize the \cst-algebra structure results of Sections \ref{kg}, \ref{ig}\ and \ref{compseq}\ for algebras of Green \ops\ 
acting between sections of vector bundles. The problem of how to define the boundary principal symbol has been addressed in \cite{RS}, 
2.3.3.1, and in \cite{G}, page 228. The crucial norm estimate (\ref{rs}) is proven in \cite{RS}, 2.3.4.4, already in this more 
general formulation.

For the K-theory applications of Sections~\ref{kt}, \ref{index}, and~\ref{superficie}, it is enough to consider the case of trivial bundles. Indeed, we are 
now going to prove that the corresponding algebras for any two choices of vector bundles are strongly Morita equivalent (this was suggested to us by Wodzicki). 
But then, strongly Morita equivalent \cst-algebras have the same K-theory. That follows, for \cst-algebras possessing countable approximate 
identities, from a theorem of Brown, Green and Rieffel \cite{BGR}. Exel \cite{E}\ proved the general case, by explicitly constructing the isomorphism 
predicted by Brown, Green and Rieffel.

Let $E$\ be a smooth rank-$k$\ vector bundle over $X$, and $F$\ be a rank-$l$\ smooth vector bundle over $\bx$, $k>0$\ and $l\geq 0$ (the case $E=0$\ belongs 
to classical index theory on closed manifolds). Let $\as\ef$\ denote the algebra of all Green \ops\
$
A:C^\infty(E)\oplus C^\infty(F)\to C^\infty(E)\oplus C^\infty(F)
$, 
and let $\ac\ef$\ denote the norm closure of $\as\ef$\ in ${\mathfrak L}(\h\ef)$, the algebra of the bounded \ops\ on 
$\h\ef=L^2(E)\oplus H^{-\frac{1}{2}}(F)$. We are going to show that $\ac\ef$\ and $\ac_{kl}$\ are strongly Morita equivalent, with $\ac_{kl}$\ denoting 
the closure of the algebra of all Green \ops\ acting between sections of the rank-$k$\ trivial bundle over $X$\ and the rank-$l$\ trivial bundle over $\bx$.

Let $\xs$\ denote the set of all Green \ops\ mapping sections of $E$\ and $F$\ to sections of the trivial bundles of rank $k$\ over $X$\ and of rank $l$\ 
over $\bx$, and let $\xc$\ denote the closure of $\xs$\ in ${\mathfrak L}(\h\ef,\h_{kl})$, $\h_{kl}=L^2(X;\co^k)\oplus H^{-\frac{1}{2}}(\bx;\co^l)$. We may 
define operator-valued inner products on $\xc$\ by $\langle A,B\rangle_{\ac_{kl}}=AB^*$\ and $\langle A,B\rangle_{\ac\ef}=A^*B$. Equipped with 
$\langle\cdot,\cdot\rangle_{\ac_{kl}}$\ and $\langle\cdot,\cdot\rangle_{\ac\ef}$, it follows immediately from the rules of \bdm's calculus that $\xc$\ becomes 
a Hilbert $\ac_{kl}$-$\ac\ef$-bimodule. To prove the strong Morita equivalence (as defined in \cite{E}, Section~5) of these two algebras, it is enough to show 
that $\xc$\ is simultaneously left and right-full. In other words, we must show that the linear span of all $A^*B$, with $A$\ and $B$\ in $\xc$, is dense in 
$\ac\ef$; and that the linear span of all $AB^*$\ is dense in $\ac_{kl}$. 

Let there be given arbitrary $C\in\as\ef$, and $\varphi$\ and $\psi\in\cix$ such that their supports are contained in a trivial open set for $E$\ 
and the intersection of their supports with $\bx$\ is contained in a trivial open set for $F$. To prove the first (the other is analogous) density statement, 
above, it clearly suffices to obtain $A$\ and $B$\ in $\xs$\ such that $\varphi C\psi=A^*B$\ (given $\rho\in\cix$, we denote also by $\rho$\ the operator 
$(f,g)\in C^\infty(E)\oplus C^\infty(F)\mapsto (\rho f, (\rho|_{\bx})g)$). This can be done by letting $B$\ have the same local expression as 
$\varphi C\psi$\ and letting $A$\ be locally given by $\tilde\psi$, for some $\tilde\psi\in\cix$\ with support slightly larger than $\supp\psi$\ and equal 
to one on a neighborhood of it. That shows that $\ac\ef$\ and $\ac_{kl}$\ are strongly Morita equivalent.

Finally, it is straightforward to prove that $\ac_{kl}$\ and $\ac$\ are strongly Morita equivalent, taking for the Hilbert bimodule the closure of
set of all Green \ops\ from $C^\infty(X,\co^k)\oplus C^\infty(\bx,\co^l)$\ to $C^\infty(X)\oplus C^\infty(\bx)$. One nontrivial step, needed only to deal 
with the case $l=0$, is to prove that every $S$\ belongs to the linear span of all $TK$, for $T$, $K$\ and $S$\ as in (\ref{boutet}).

For many purposes, one might therefore assume that the bundle
over the boundary is zero. In order to reach the statements of
Corollaries 12 and 20, for example, the reader may, from now on,
focus on the upper left corner of the matrix in \eqref{boutet} and
omit all statements referring to the other entries. 

%
%
%
%
%
\section{The kernel of $\gam$}\label{kg}

Let us first note that $\is$\ is contained in the kernel of $\gamma$, as implied by the facts that $G$, $K$, $T$\ and $S$\ in 
(\ref{dis}) are of negative order, and that the principal symbol of $\varphi P\psi$\ vanishes over the boundary.

Another easy observation is that $\ic$\ contains the compact ideal $\kc$. That follows because the set of all \ops\ with smooth kernel is dense in $\kc$\
and any such \op\ can be written as in (\ref{dis}), with $P=0$\ and $G$, $K$, $T$\ and $S$\ of order $-{\infty}$.
 
\begin{lem}\label{l1} Every $A\in\as$\ such that $\gamma(A)=0$\ belongs to $\ic$. \end{lem}

\pf Let $A\in\as$\ be given, with $\gamma(A)=0$. It follows from the definition that $A_{jk}\in\is_{jk}$, if $j\neq 1$\ or $k\neq 1$. 
If $A_{11}=\pp+G$, then $\gamma(\pp)(\cvb)=-\gamma(G)(\cvb)$\ for every $(\cvb)\in\sbx$. In Section~\ref{intr}, we remarked that $\gamma(G)(\cvb)\in\kc$\ for 
every $(\cvb)\in\sbx$\ and that $p(D)_{_{+}}\mapsto p$\ defines an isomorphism of $\tc/\kc$\ with $C(\re\cup\{\infty\})$. From that it follows that the 
principal symbol of $P$, $p_0\in\cisx$, vanishes over the boundary $\bx$\, and, hence, that $\gamma(\pp)$\ and $\gamma(G)$\ both vanish. Then, by (\ref{rs}), 
$G$\ is compact ($\sigma(G)=0$, by definition) and it suffices to show that $\pp$\ is in $\ic_{11}$. 

Let $\rho\in\cix$\ be a {\em boundary defining function}, i.e., $\rho$\ is positive on $\ix$\ and vanishes with nonzero derivative on $\bx$. Let $Q$\ be a 
zero-order classical \psd\ \op\ on $X$\ with principal symbol $q\in\cisx$\ defined by $p_0=\rho q$. We have $\pp\equiv\rho Q_{_{+}}$, modulo a compact \op, 
since \psd\ \ops\ of negative order on compact manifolds are compact. Moreover, $Q$\ has the \tp, but we will not need that. Now let 
$\varphi_k\in\cciix$\ be equal to one on some sequence of compacts exhausting $\ix$. Since $\varphi_k\rho$\ converges to $\rho$\ uniformly on $X$, we have 
that $\varphi_k\rho Q_{_{+}}$\ converges to $\rho Q_{_{+}}$\ in ${\mathfrak L}(L^2(X))$. 

So, it suffices to prove that $\varphi Q_{_{+}}\in\ic_{11}$\ if $\varphi$\ belongs to $\cciix.$\ Let $\psi\in\cciix$\ be equal to one in a neighborhood of 
$\varphi$. Then we have: $\varphi Q_{_{+}}=\psi[\varphi,Q_{_{+}}]+\psi Q_{_{+}}\varphi$. The commutator $[\varphi,Q_{_{+}}]$\ is compact, since it is equal 
to $(\varphi Q-Q\varphi)_{_{+}}$\ and $(\varphi Q-Q\varphi)$\ has negative order. That proves the lemma, since $\psi Q_{_{+}}\varphi=\psi Q\varphi$. \cqd

Thus, we have $\is\subset\ker\gamma\subset\ic$, and, hence, $\overline{\ker\gamma}=\ic$. It is obvious that $\overline{\ker\gamma}\subseteq\ker\gam$. We prove 
next that $\ker\gam\subseteq\ic$. Let $A\in\ker\gam$\ be given. If $j\neq1$\ or $k\neq1$, we then have:
\begin{equation}
\label{fullaint}
\inf_{A^\prime\in\ic}||A_{jk}+A_{jk}^\prime||\leq\inf_{C\in\kc}||A_{jk}+C||=\max\{||\sig_{jk}(A)||,||\gam_{jk}(A)||\}=||\gam_{jk}(A)||,
\end{equation}
since, by definition, $\sig_{jk}=0$\ unless $j=k=1$. We have used (\ref{rsc}) for the matrix that has $A_{jk}$\ in its 
$(j,k)$-entry and zero in the others. It follows from (\ref{fullaint}) that $\ker\gam_{jk}\subseteq\ic_{jk}$\ if $j\neq1$\ or $k\neq 1$.

To show that $\ker\gam_{11}\subseteq\ic_{11}$\ (and hence prove that $\ker\gam=\ic$), it is enough to prove Lemma~\ref{crucial}, below. Indeed, the left-hand side
of (\ref{trago}) is by definition greater than or equal to $\inf\{||\pp+G+A^{\prime}_{11}||;A^\prime\in\ic\}$. Moreover, since both sides of the inequality 
depend continuously on $A=\pp+G$, proving (\ref{trago}) will imply that it holds also, for the extension $\gam$, if $\pp+G$\ is replaced by a general 
$A\in\ac_{11}$.

For each $s>0$, let us define the unitary \op\ $\kappa_s$\ on $L^2(\rp)$\ by $(\kappa_s\psi)(t)=\sqrt{s}\psi(s t)$. It is straightforward to check that
\begin{equation}
\label{Fourier}
\kappa_s^{-1}b(D)_{_{+}}\kappa_s\,=\,b(sD)_{_{+}},\ \text{for all}\ s>0,
\end{equation}
with $b(D)_{_{+}}$\ denoting the \op\ on $L^2(\rp)$\ obtained by truncating the Fourier multiplier
$b(D)$, for any measurable function $b$\ on $\re$. 

\begin{lem}
\label{crucial}
There exists a positive constant $c$, determined only by $X$, such that, for every \psd\ \op\ with the \tp\ on $X$, and for
every singular Green \op\ $G$, both polyhomogeneous and of order zero, we have:
\begin{equation}
\label{trago}
\inf_{\varphi,\psi,Q}||\pp+G+\varphi Q\psi||\leq c||\gamma(\pp+G)||,
\end{equation}
where the infimum is taken over all $\varphi$\ and $\psi$\ in $\cciix$\ and all zero-order classical \psd\ \ops\ $Q$\ on $X$.
\end{lem}
\pf The closure of set of all such $\varphi Q\psi$\ contains the compact \ops\ (it is enough to take $Q$'s of order $-\infty$). Moreover, (\ref{rs}) 
applied to the matrix that has $G$\ in the upper left corner and zero in the other entries gives $\inf_{C\in\kc}||G+C||\,=\,||\gamma(G)||$. 
This implies (\ref{trago}), with $c=1$, for the case $\pp=0$.

Next we prove (\ref{trago}) for the case $G=0$. We have already noticed that the left-hand side of (\ref{trago}) is not smaller than the infimum of 
$\{||\pp+G+A^\prime_{11}||;A^\prime\in\ic\}$. 
They are, in fact, equal, since any singular Green \op\ of negative order on $X$\ is compact, and the compacts are contained in the closure of the set of 
all $\varphi Q\psi$. Let $\{\varphi_1,\cdots,\varphi_d\}$\ be a partition of unity on $X$\ such that, whenever the supports of 
$\varphi_i$\ and $\varphi_j$\ intersect, their union is contained in the domain of a chart $\chi_{ij}:U_{ij}\to\tilde U_{ij}$. If $\supp\varphi_i$\ and 
$\supp\varphi_j$\ do not intersect, then $\varphi_i\pp\varphi_j$\ is regularizing, hence compact. Since ${\mathfrak K}\subseteq\ic_{11},$\ we get:
\begin{equation}
\label{local}
\inf_{A^\prime\in\ic_{11}}||\pp+A^\prime||\leq\sum_{\supp\,\varphi_{i}\cap\supp\,\varphi_{j}\neq\emptyset}\inf_{A^\prime\in\ic_{11}}
||\varphi_i\pp\varphi_j+A^\prime||.
\end{equation}
If $U_{ij}$\ does not intersect the boundary, then $\varphi_i\pp\varphi_j=\varphi_iP\varphi_j\in\is_{11}$. For each $(i,j)$\ in the sum, above, we may 
therefore suppose that $U_{ij}$\ intersects the boundary, given by $x_n=0$, and denote by $p(x,\xi)$\ the local symbol of $P$\ for that chart. Let $P_{ij}$\ 
denote the pullback by $\chi_{ij}$\ of the \psd\ \op\ $\tilde P_{ij}$\ on $\rnp$\ of amplitude 
\begin{equation}
\label{qij}
q_{ij}(x,y,\xi)=p(x^\prime,0,\xi)\varphi_{i}(\chi_{ij}^{-1}(x))\varphi_{j}(\chi_{ij}^{-1}(y)).
\end{equation}
$P_{ij}$\ is a classical \psd\ \op\ with the \tp\ on $X$, such that $P_{ij_{+}}-\varphi_i\pp\varphi_j$\ is in the kernel of $\gamma$. 
It follows from Lemma~\ref{l1}\ that $P_{ij_{+}}-\varphi_i\pp\varphi_j$\ belongs to $\ic_{11}$\ and, hence, that 
\begin{equation}
\label{lcl}
\inf_{A^\prime\in\ic_{11}}||\varphi_i\pp\varphi_j+A^\prime||=\inf_{A^\prime\in\ic_{11}}||P_{ij_{+}}+A^\prime||\leq\inf_{C\in\kc}||P_{ij_{+}}+C||.
\end{equation}

Let us denote by $P^\prime_{ij}$\ the \psd\ \op\ on $\rn$\ with amplitude defined by the same formula (\ref{qij}) as $\tilde P_{ij}$, simply assuming, as we
may, that the $\varphi$'s and $\chi$'s are restrictions of functions and charts on the neighboring manifold $\Omega$. The classical estimate for the norm, 
modulo compacts, of a \psd\ 
\op\ in terms of the supremum-norm of its principal symbol (\cite{KN}, Theorem A.4; or \cite{Hp}, Theorem 3.3) implies the existence of compact \ops\ 
$C^\prime_{ij}$\ on $L^2(\rn)$, such that $||P^\prime_{ij}+C^\prime_{ij}||$\ is bounded by two times the supremum-norm of the principal symbol of 
$P^\prime_{ij}$; which is bounded by the supremum on the right-hand side of (\ref{nev}), below, since $|\varphi_k|\leq 1$\ for all $k$. If $\tilde C_{ij}$\
denotes the compact \op\ on $L^2(\rnp)$\ obtained from $C^\prime_{ij}$\ by truncation, it is obvious that $||\tilde P_{ij_{+}}+\tilde C_{ij}||$\ is bounded by 
$||P^\prime_{ij}+C^\prime_{ij}||$. This gives: 
\begin{equation}
\label{nev} 
||\tilde P_{ij_{+}}+\tilde C_{ij}||\leq 2\sup\{|p_0(x^\prime,0,\xi)|;\,(x^\prime,0)\in\tilde U_{ij},\,\xi\neq 0\},
\end{equation}
where $p_0(x,\xi)$\ (smooth for $\xi\neq 0$) denotes the zero-order homogeneous principal part of $p(x,\xi)$. 

This is the most delicate point of this proof: 
For each $(\cvb)$\ with $\xi^\prime\neq 0$, $p_0(x^\prime,0,\xi^\prime,\cdot)$\ belongs to ${\mathcal H}_{0}$, and 
\begin{equation}
\label{estimate}
\sup_{\xi_{n}\in\re}|p_0(x^\prime,0,\xi^\prime,\xi_n)|=||p_0(x^\prime,0,\xi^\prime,D_n)_{_{+}}||_{{\mathfrak L}(L^2(\rp))}.
\end{equation}
That again follows from the isomorphism $\tc/\kc\simeq C(\re\cup\{\infty\})$, or, in a more classical language, from Lemma~3.1.5 of \cite{G}. Since
\[
\sup_{\xi^\prime\neq 0}\sup_{\xi_n\in\re}|p_0(x^\prime,0,\xi^\prime,\xi_n)|\,=\,\sup_{\xi\neq 0}|p_0(x^\prime,0,\xi)|,
\]
the right-hand side of (\ref{nev}) equals $2\sup_{\xi^\prime\neq 0}||p_0(x^\prime,0,\xi^\prime,D_n)_{_{+}}||$. It follows from the homogeneity of $p_0$\
and (\ref{Fourier}) that $p_0(x^\prime,0,s\xi^\prime,D_n)_{_{+}}=\kappa_sp_0(x^\prime,0,\xi^\prime,D_n)_{_{+}}\kappa_s^{-1}$, for all $s>0$, and hence
$||p_0(x^\prime,0,\xi^\prime,D_n)_{_{+}}||$ is independent of $|\xi^\prime|$. We then get:
\begin{equation}
\label{verao}
||\tilde P_{ij_{+}}+\tilde C_{ij}||\leq 2||\gamma(\pp)||.
\end{equation}

Let $C_{ij}\in\kc$\ denote the pullback of (a restriction of) $\tilde C_{ij}$\ by the chart $\chi_{ij}$ ($C_{ij}$\ vanishes on functions whose 
support does not intersect the closure of $U_{ij}$). There is a constant $c_1$, depending only on our choice of norm on $L^2(X)$, such that
\begin{equation}
\label{nv}
||P_{ij_{+}}+C_{ij}|| \leq c_1
||\tilde P_{ij_{+}}+\tilde C_{ij}|| .
\end{equation}
Estimates (\ref{local}), (\ref{lcl}), (\ref{verao}) and (\ref{nv}) imply (\ref{trago}), for the case $G=0$, with $c=2c_1d^2$.

To treat the case when both $\pp$\ and $G$\ are nontrivial, we use that $\gamma(G)(\cvb)$\ is compact for each $(\cvb)\in\sbx$, and again the
isomorphism $\tc/\kc\simeq C(\so)$, to get:
\[
||\gamma(\pp)(\cvb)||=\inf_{C\in\kc}||\gamma(\pp)(\cvb)+C||\leq||\gamma(\pp)(\cvb)+\gamma(G)(\cvb)||.
\]
Taking the supremum on both sides of this inequality, we see that $||\gamma(\pp)||\leq||\gamma(\pp+G)||$, and, hence, $||\gamma(G)||\leq 2||\gamma(\pp+G)||$.
Since 
\[
\inf_{\varphi,\psi,Q}||\pp+G+\varphi Q\psi||\,\leq\,\inf_{\varphi,\psi,Q}||\pp+\varphi Q\psi||+\inf_{\varphi,\psi,Q}||G+\varphi Q\psi||,
\]
the proof is complete, with $c=2(1+c_1d^2)$. \cqd

\begin{cor}\label{ami}
$\gam$\ induces a \cst-algebra isomorphism between 
$\ac/\ic$\ and $\text{{\em Im}}\,\gam$. Equality, with $c=1$, therefore holds in (\ref{trago}).
\end{cor}

For the sake of this argument, let $\kc_X$\ and $\kc_\oplus$\ denote the ideals of compact \ops\ on $L^2(X)$\ and on $\h$, respectively. There is an obvious 
injection of $\ic_{11}/\kc_X$\ into $\ic/\kc_\oplus$; which is also surjective because, by definition, all entries of the matrix in (\ref{dis}) are compact, 
except possibly the upper left one. The statement about $\ic/\kc$\ in the theorem, below, follows therefore from the estimate
\[
\inf_{C\in\kc}||\varphi P\psi+C||\ =\ \sup_{S^{*}\ix}|\varphi\psi\sigma(P)|,
\]
which again follows from the classical estimate quoted between (\ref{lcl}) and (\ref{nev}) (see also \cite{JK}, Section 2). We have proven:

\begin{thm}\label{tker}\ The kernel of $\gam$\ is equal to $\ic$. Moreover, $\ic$\ contains the compact ideal $\kc$\ and $\ic/\kc$\ is isomorphic to $\cosx$, 
with isomorphism induced by the principal symbol. \end{thm}
%
%
%
%
%
%
%
%
%
%
\section{The image of $\gam$}\label{ig}

It is of central importance for the computation of the K-theory of \bdm's algebra that the upper left corner of the image of $\gam$, $\img_{11}$, 
is equal to the Banach-space direct sum of a subalgebra isomorphic to $C(\bx)$\ with an ideal with vanishing K-theory. This description is the essential result of this 
section; what we state about $\img_{jk}$, for $j\neq1$\ or $k\neq1$, is already contained in \cite{RS}, 2.3.4.4, Corollary~2. 

Every smooth function on $\sbx$\ is the principal symbol of a zero-order classical \psd\ \op\ on the closed manifold $\bx$. 
$\cisbx$\ is therefore contained in $\text{Im}\,\gamma_{22}$\ and, hence, $\img_{22}=C(\sbx)$. The fact that $\img_{21}=C(\sbx)\otimes L^2(\rp)^{*}$\ 
(so, $\gam_{21}$\ is surjective, if we define $\gam$\ taking values in the \cst-algebra in (\ref{gama})) follows from the next lemma, 
by a partition-of-unity argument. 

Given $\varphi\in L^2(\rp)$, we will denote by $\langle\varphi|$\ the linear functional $\psi\mapsto\int\varphi\psi$; and by $|\varphi\rangle$\ the linear
map, from $\co$\ to $L^2(\rp)$, of multiplication by $\varphi$.

%
%
%
\begin{lem}\label{trivial}
Let $V\subset\bx$\ be an open set whose closure is contained in the domain of a chart of $\bx$. The space $C_0(S^{*}V,L^2(\rp)^{*})$\ of all continuous 
functions, from $\overline{S^{*}V}$\ to $L^2(\rp)^{*}$, which vanish on every unit covector over the boundary of $V$\ is contained in the image of $\gam_{21}$.
\end{lem}

\pf Let $\chi:U\to\tilde{U}\subseteq\re^{n-1}$\ be a chart of $\bx$\ such that $\bar{V}\subset U$.
Let there be given $p$\ in $C_{c}^{\infty}(S^{*}V)$\ (i.e., $p\in\cisbx$\ and its support is contained in $S^{*}V$) and $\varphi$\ in 
$C_c^{\infty}(\re_+)$.\ Let us denote by $\tilde p(\cvb)$, $(\cvb)\in\tilde{U}\times\re^{n-1}$, the local expression (smooth for $\xi^\prime\neq 0$) of the 
zero-degree homogeneous extension of $p$\ to the cotangent bundle of $\bx$. We now choose an {\em excision function}\ $\omega$\ (i.e., 
$\omega\in C^\infty(\re)$\ vanishes on a neighborhood of the origin and $\omega(t)\equiv 1$\ for sufficiently large $t$) and define 
$t(\cvb,\xi_n)=2\pi\omega(|\xi^\prime|)\tilde p(\cvb)\check\varphi(\xi_n/|\xi^\prime|)$, with $\,\check{}\,$\ denoting the inverse Fourier transform, and
$|\xi^\prime|$\ the euclidean norm of $\xi^\prime\in\re^{n-1}$. 

Obviously, $t$\ is smooth. One can also check that $t$\ is a trace symbol of order and class zero (using \cite{G}, (1.2.19) and (2.3.25), for example), 
defining therefore a trace \op\ $\tilde T:C_c^{\infty}(\rnp)\to C_c^{\infty}(\re^{n-1})$. Moreover, $t$\ is polyhomogeneous and its homogeneous principal 
part is given by $t_0(\cvb,\xi_n)=2\pi\tilde p(\cvb)\check\varphi(\xi_n/|\xi^\prime|)$. We then get (using \cite{G}, (2.4.5) and (2.3.25), for example): 
\begin{equation}
\label{tracesymbol}
t_0(\cvb,D_n)\,=\,|\xi^\prime|\,\tilde p(\cvb)\langle\varphi(|\xi^\prime|\cdot)|.
\end{equation}

Let $D$\ denote the linear span of all $p\otimes\langle\varphi|$, with $p\in C_c^{\infty}(S^{*}V)$, $\varphi\in C_c^{\infty}(\re_+)$. Since 
$C_c^{\infty}(\re_+)$\ is dense in $L^2(\rp)$, $D$\ is dense in $C_0(S^{*}V,L^2(\rp)^{*})$. The assignment
\begin{equation}
\label{iota}
D\ni f(\cvb)\mapsto \sqrt{|\xi^\prime|}f(\cvb)\circ\kappa_{|\xi^\prime|}
\end{equation}
($\kappa_{(\cdot)}$\ as defined before Lemma~\ref{crucial}), induces an isomorphism $\iota$\ of $C_0(S^{*}V,L^2(\rp)^{*})$\ onto itself, since 
it is equal to the multiplication of an isometry by the function $\sqrt{|\xi^\prime|}$, which is smooth, bounded, and bounded away from zero on $S^{*}V$. 

Let $\rho\in C_c^\infty(\tilde V)$, $\tilde V=\chi(V)$, be equal to one in a neighborhood of the $x^\prime$-support of $\tilde  p$. Then, extend $\chi$\ to a 
chart $\chi_1$\ of $X$, and $\rho$\ to a cutoff $\rho_1$\ with support in the image of $\chi_1$.
The pullback $T:\cix\to\cibx$\ of $\rho\tilde T\rho_1$\ by $\chi_1$\ is a trace \op\ and  $\gamma(T)=\iota(p\otimes\langle\varphi|)$. Hence, $\img_{21}$\ 
contains $\iota(D)$, which is dense in $C_0(S^{*}V,L^2(\rp)^{*})$. On the other hand, $\img_{21}$\ is closed, since $\gam$\ is a \cst-algebra homomorphism. 
\cqd

The surjectivity of $\gam_{12}$\ can be proven analogously to that of $\gam_{21}$. The main difference is that the homogeneous extension of $p$, which in the proof of 
Lemma~\ref{trivial}\ was of degree zero, now has to be taken of degree $(-1)$. One should then define the local Poisson symbol by 
$k(\cvb,\xi_n)=2\pi\omega(|\xi^\prime|)\tilde p(\cvb)\check\psi(\xi_n/|\xi^\prime|),$\ with $\psi(t)=\varphi(-t)$; and replace $\langle\cdot|$\ by
$|\cdot\rangle$\ in (\ref{tracesymbol}), and $f(\cvb)\circ\kappa_{|\xi^\prime|}$\ by $\kappa_{|\xi^\prime|}\circ f(\cvb)$\ on the right-hand side of 
(\ref{iota}).

We must now describe the upper-left corner of $\img$. As a first step, the following lemma will show that $\img_{11}$\ contains $C(\sbx)\otimes\tc_0$ 
($\tc_0=\ker\lambda$, as defined after (\ref{gama})\,). Before proving it, let us show that the \cst-algebra (let us call it $\tc_1$, for the moment)
generated by all $\varphi(D)_{_{+}}$, $\varphi\in{\mathcal S}(\re)$, is equal to $\tc_0$. Since every generator of $\tc_1$\ clearly belongs to $\tc_0$, one
gets at once that $\tc_1\subseteq\tc_0$. 

Given $\varphi$\ in the space $C_0(\re)$\ of all continuous functions on $\re$\ that vanish at infinity, $U^{-1}F\varphi(D)_{_{+}}F^{-1}U$\ is equal to
the Toeplitz \op\ on $\so$\ of symbol $\varphi(\frac{iz-i}{z+1})$\ (recall that $U$\ was defined in (\ref{unitary}) and $F$\ denotes the Fourier transform). 
The isomorphism $\tc/\kc\simeq C(\so)$\ then implies that $\tc_0$\ is equal to the set of all $\varphi(D)_{_{+}}+C$, with $\varphi\in C_0(\re)$\ and $C$\ 
compact. It follows from the estimate $||\varphi(D)_{_{+}}||=\sup|\varphi|$\ that every $\varphi(D)_{_{+}}$, with $\varphi$\ in $C_0(\re)$, is the norm limit 
of a sequence $\varphi_k(D)_{_{+}}$, $\varphi_k\in{\mathcal S}(\re)$. This shows that $\tc_1$\ is equal to the \cst-algebra generated by all 
$\varphi(D)_{_{+}}$, $\varphi\in C_0(\re)$. To prove that $\tc_0\subseteq\tc_1$, it is therefore enough to show that $\tc_1$\
contains $\kc$. This follows from the fact (\cite{Dou}, Proposition~7.12) that  $\kc$\ is equal to the commutator ideal of $\tc$, which is  
equal to the commutator ideal of $\tc_1$. Indeed, the two ideals are equal to the closed linear span of all products $T_1\cdots T_k$, where at least one
(possibly more) of the $T_j$'s is of the form $[\varphi(D)_{_{+}},\psi(D)_{_{+}}]$, $\varphi$\ and $\psi$\ in $C_0(\re)$, and the others are of the form
$\varphi(D)_{_{+}}$, $\varphi$\ in $C_0(\re)$. 

\begin{lem}
\label{tz}
Let $V$\ be as in Lemma~\ref{trivial}. Then $C_0(S^{*}V,\tc_0)$\ is contained in $\img_{11}$.
\end{lem}
\pf 
Given $\varphi\in{\mathcal S}(\re)$\ and $p\in C_c^{\infty}(S^{*}V)$, let $P$\ denote the pullback by $\chi_1$\ of the \psd\ \op\ on $\rnp$\ of amplitude 
\begin{equation}
\label{onto11}
q(x,y,\xi)=\omega(|\xi^\prime|)\tilde p(\cvb)\varphi(\xi_n/\sqrt{1+|\xi^\prime|^2})\rho_1(x)\rho_1(y),
\end{equation}
with $\chi$, $\chi_1$, $\omega$, $\tilde p$\ and $\rho_1$\ as in the proof of Lemma~\ref{trivial}. It follows from the proof of Lemma~5.3.1 in 
\cite{schsch}\ that $q$\ is a classical amplitude, and that the corresponding homogeneous principal symbol $q_0(x,\xi)$\ satisfies, for 
$x^\prime\in\tilde V$\ and $\xi\neq 0$, $q_0(x^\prime,0,\xi)=\tilde p(\cvb)\varphi(\xi_n/|\xi^\prime|)$. 
Moreover, $P$\ has the transmission property (\cite{RS}, 2.2.2.1). We then get from (\ref{Fourier}):
\begin{equation}
\label{generator}
\gamma(\pp)(\cvb)\,=\,\tilde p(\cvb)\cdot[\kappa_{|\xi^\prime|}\circ\varphi(D)_{_{+}}\circ\kappa_{|\xi^\prime|}^{-1}].
\end{equation}

It follows from the fact that $\tc_0$\ is generated by all $\varphi(D)_{_{+}}$, $\varphi\in{\mathcal S}(\re)$, that $C_0(S^{*}V,\tc_0)$\ is equal to the 
\cst-algebra generated by all $\tc_0$-valued functions on $S^{*}V$\ as those at the right-hand side of (\ref{generator}), with $\varphi\in{\mathcal S}(\re)$\ 
and $p\in C_c^{\infty}(S^{*}V)$. \cqd 

Now we show that $\gam$\ is not surjective. We will regard $C(\bx)$\ as a subset of $C(\sbx)$, and $C(\sbx)$\ as a subset of $C(\sbx)\otimes\tc$\ 
(recall that $\tc$\ contains the identity \op\ $I$\ on $L^2(\rp)$\,).

\begin{lem}
\label{notonto}
{\em $\img_{11}\cap C(\sbx)\ =\ C(\bx). $}
\end{lem}
\pf
Operators of multiplication by functions in $\cix$\ are the simplest examples of \psd\ \ops\ with the \tp\ on $X$. The boundary principal symbol of such an 
\op\ is 
equal to the restriction to the boundary of the multiplier (times $I$). This shows that $\cibx$, and hence also $C(\bx)$, are contained in $\img_{11}$.

Let $f$\ belong to $\img_{11}\cap C(\sbx)$. For every $\delta>0$, there exist a \psd\ \op\ with the \tp\ $P$ and a singular Green \op\ $G$ such that 
$||\gamma(\pp)(\cvb)+\gamma(G)(\cvb)-f(\cvb)||<\delta$, for all $(\cvb)\in\sbx$\ (we have used that the set of all $\pp+G$\ is dense in $\ac_{11}$). Since 
$\gamma(G)(\cvb)$\ is compact for every $(\cvb)$, we get:
\begin{equation}
\label{delta}
\inf_{C\in\kc}||\gamma(\pp)(\cvb)-f(\cvb)I+C||\,<\,\delta.
\end{equation}

Using once more the isomorphism $\tc/\kc\simeq C(\re\cup\{\infty\})$, we get 
\[
\inf_{C\in\kc}||p(D)_{_{+}}+C||\,=\,\sup|p|,
\] 
for any $p\in{\mathcal H}_0$. Let $p_0$\ denote the principal symbol of $P$, regarded as a zero-degree homogeneous function on the cotangent bundle,
smooth except at the zero section. The left-hand side of (\ref{delta}) is then equal to 
\[
\sup_{\xi_{n}\in\re}|p_0(x^\prime,0,\xi^\prime,\xi_n)-f(\cvb)|=\sup_{\xi_{n}\neq 0}|p_0(x^\prime,0,\frac{\xi^\prime}{|\xi_{n}|},\pm 1)-f(\cvb)|,
\]
which is greater than or equal to $|p_0(x^\prime,0,0,+1)-f(\cvb)|$. This makes sense in view of our choice of a normal coordinate $x_n$\ after 
(\ref{pbs}). Note also that, for covectors $\xi^\prime\neq 0$, $p_0(x^\prime,0,\xi^\prime,\cdot)\in{\mathcal H}_0$.

We have proven that, for all $\delta>0$, there exists a $g\in C(\bx)$\ such that $\sup_{\sbx}|f-g|<\delta$. 
Hence, $f\in C(\bx)$. \cqd

The previous argument also shows that $p_0(x^\prime,0,0,+1)=p_0(x^\prime,0,0,-1)$. This is part of \bdm's transmission condition for classical 
\ops, see the comments before the statement of Theorem~\ref{stil}\ for more details. 

Lemmas \ref{tz}\ and \ref{notonto}\ show that $\img_{11}\subseteq (C(\sbx)\otimes\tc_0)\oplus C(\bx)$. To prove the reverse inclusion, 
let us consider the \cst-algebra homomorphism $1\otimes\lambda$\ from $C(\sbx)\otimes\tc$\ to itself that maps $f\otimes p(D)_{_{+}}$\ to 
$p(\infty)f\otimes I$, $p\in{\mathcal H}_0$. If $F$\ belongs to $\img_{11}$, then $F-(1\otimes\lambda)(F)$\ belongs to 
$C(\sbx)\otimes\tc_0$. By Lemma~\ref{tz}, $F-(1\otimes\lambda)(F)$, and therefore also $(1\otimes\lambda)(F)$, belong to $\img_{11}$. By Lemma~\ref{notonto}, 
$(1\otimes\lambda)(F)$\ is in $C(\bx)$. 

This proves the characterization of $\img$\ promised after (\ref{gama}):
\begin{equation}
\label{imofga}
\img=\left(C(\sbx)\otimes\left(\begin{array}{cc}\tc_0&L^2(\rp)\\L^2(\rp)^{*}&\co\end{array}\right)\right)\oplus
\left(C(\bx)\otimes\left(\begin{array}{cc}\co&0\\0&0 \end{array}\right)\right). 
\end{equation} 

Before stating the main result of this Section, however, a few more definitions and comments are needed. Let us denote by $\wc$\ the \cst-algebra of all the 
bounded \ops\ $A=(\!(A_{jk})\!)_{j,k=1,2}$\ on $L^2(\rp)\oplus\co$\ such that $A_{11}\in\tc$, and by $\wc_0$\ the set of all $A\in\wc$\ such that 
$A_{11}\in\tc_0$\ ($\wc$\ is the closure of the set of all zero-order {\em Wiener-Hopf \ops}, in \bdm's terminology \cite{B2}). 
$\wc_0$\ is clearly the kernel of the linear functional $\Lambda(A)=\lambda(A_{11})$; which is actually a homomorphism, since 
$\Lambda(AB)=\lambda(A_{11}B_{11}+A_{12}B_{21})$\ and $A_{12}B_{21}$\ is compact. 

Temporarily denoting by $\kc_{\rp}$\ the ideal of compact \ops\ on $L^2(\rp)$\ and by $\kc_{\oplus}$\ that on $L^2(\rp)\oplus\co$, it is straightforward 
to check that the mapping that sends the class of $A$\ to the class of the matrix that has $A$\ in the upper-left corner and zero elsewhere is an 
isomorphism between $\tc/\kc_{\rp}$\ and $\wc/\kc_{\oplus}$. The isomorphism $\tc/\kc_{\rp}\simeq C(\so)$\ defined in the introduction therefore induces an 
isomorphism from $\wc/\kc_{\oplus}$\ to $C(\so)$\ that sends the class of the identity $I$\ to the constant function $1$. The homomorphism $\Lambda$\ 
corresponds, then, to the evaluation at $(-1)$\ on $C(\so)$\ composed with the projection of $\wc$\ onto $\wc/\kc_{\oplus}$.

We have proven:

\begin{thm}
\label{imagem} The \cst-algebra $\wc$\ contains the compact ideal, and $\wc/\kc$\ and $C(\so)$\ are isomorphic (as unital \cst-algebras).  
Moreover, the image of $\gam$\ is isomorphic, as a Banach space, to $C(\bx)\oplus (C(\sbx)\otimes\wc_0)$, with $\wc_0$\ denoting the kernel of the homomorphism
$\Lambda:\wc\to\co$\ induced by the evaluation at $(-1)$\ on $C(\so)$.
\end{thm}
%
%
%
%
%
%
%
%
\section{K-Theory}\label{kt}

By a well-known consequence of Bott periodicity  (\cite{Bl}, Theorem 9.3.1), to any short exact sequence of \cst-algebras $0\to J\to A\to A/J\to 0$, 
one may associate a cyclic six-term exact sequence of abelian groups
\begin{equation}
\label{sixterm}
\begin{array}{ccccc}
K_0(J)&\longrightarrow&K_0(A)&\longrightarrow&K_0(A/J)\\\big\uparrow&&&&\big\downarrow\\K_1(A/J)&\longleftarrow&K_1(A)&\longleftarrow& K_1(J)
\end{array},
\end{equation}
where the horizontal arrows are functorially induced homomorphisms, the arrow connecting $K_1$\ to $K_0$\ is called the {\em index mapping}, and the other 
connecting mapping is called the {\em exponential mapping}. If $J$\ is the compact ideal, then $K_1(J)=0$, $K_0(J)\simeq\z$, and the index mapping 
is the Fredholm index (\cite{Bl}, 8.3.2). 

\begin{lem}
\label{ktht}
$K_i(C(\sbx)\otimes\wc_0)=0$, $i=0,1$.
\end{lem}
\pf The six-term exact sequence we get from the short exact sequence given by the isomorphism of our Theorem~\ref{imagem}, 
$0\to\kc\to\wc\to C(\so)\to 0$, is 
\begin{equation}
\label{T6}
\begin{array}{ccccc}
\z&\longrightarrow&K_0(\wc)&\longrightarrow&\z\\\big\uparrow&&&&\big\downarrow\\\z&\longleftarrow&K_1(\wc)&\longleftarrow& 0.
\end{array}
\end{equation}
Because there exists a Toeplitz \op\ of Fredholm index one, there exists also an \op\ in $\wc$\ of index one (take a matrix with 1 in the lower right corner
and 0 outside the diagonal). Hence, the index mapping in (\ref{T6}) is surjective. This gives 
$
K_0(\wc)=[I]\cdot\z\ \text{and}\ K_1(\wc)=0.
$
It then follows from the six-term exact sequence associated to  
\[
0\longrightarrow\wc_0\longrightarrow\wc\smash{\mathop{\longrightarrow}\limits^\Lambda}\co\longrightarrow 0 
\] 
that $K_{0}(\wc_0)=K_{1}(\wc_0)=0$. By K\"unneth's theorem for tensor products \cite{Bl,WO}, we prove our claim. \cqd

Denoting by the same letters functions on $X$\ or $\bx$\ and the multiplication \ops\ they define, let $b$\ denote the (unital) \cst-algebra homomorphism
\begin{equation}
\label{bdmp}
\begin{array}{rcl}
b:C(\bx)&\longrightarrow&\img\\
g&\longmapsto &\gam\left(\left(
\begin{array}{cc}f&0\\0&g\end{array}
\right)\right),
\end{array}
\end{equation}
where $f$\ denotes a function in $C(X)$\ whose restriction to the boundary equals $g$.

\begin{cor}
\label{kthbx}
The homomorphisms $b_*:K_{i}(C(\bx))\to K_{i}(\text{{\em Im}}\,\gam)$, $i=0,1$, induced by $b$\ are isomorphisms.
\end{cor}

\pf  The proof of Theorem~\ref{imagem}\ also shows that $\img\simeq(C(\sbx)\otimes\wc_0)\oplus
\text{Im}\,b$ (this follows from  \eqref{imofga}).
With respect to this Banach-space direct-sum decomposition, 
the product on $\img$\ is given by
\[
(F\oplus b(f))(G\oplus b(g))=(FG+fG+gF)\oplus b(fg).
\]
In particular, ${\mathfrak R}=C(\sbx)\otimes\wc_0$\ is an ideal of $\img$, and $\img/{\mathfrak R}\simeq\text{Im}\,b$. The six-term exact sequence associated to
$0\to{\mathfrak R}\to\img\to\text{Im}\,b\to 0$, together with Lemma~\ref{ktht}, imply that $\pi_*:K_i(\img)\to K_i(\text{Im}\,b)$, $i=0,1$, are isomorphisms, with
$\pi$\ denoting the canonical projection of $\img\simeq{\mathfrak R}\oplus\text{Im}\,b$\ onto $\text{Im}\,b$. 

Since $b$\ is injective, $\pi\circ b:C(\bx)\to\text{Im}\,b$\ is a \cst-algebra isomorphism. We then get $b_*$\ equal to 
the composition of group isomorphisms $\pi_*^{-1}\circ(\pi\circ b)_*$. \cqd

\begin{pro}
\label{tritan}
If $X$\ is connected and $\bx$\ is not empty, then there exists a nonvanishing section of the cotangent bundle of $X$.
\end{pro}
\pf Let $\Sigma$\ be a section of the cotangent bundle of $\Omega$\ with a finite number of zeros. Given $x_0\in X$, a zero of $\Sigma$, let $c:[0,1]\to\Omega$\ 
be a smooth curve with $c(0)=x_0$\ and $c(1)\not\in X$, such that $c(t)$\ is not a zero of $\Sigma$\ if $t\neq 0$, and $c^\prime(t)\neq 0$\ for all $t$.
Let $V$\ be a vector field on $\Omega$\ such that $V(c(t))=c^\prime(t)$, $t\in[0,1]$, and $V\equiv 0$\ outside a neighborhood of the image of $c$\ which 
intersects the set of zeros of $\Sigma$\ only at $x_0$. Let $\Phi_s$, $s\in\re$, denote the flow of $V$, and write $\tilde\Sigma$\ for the pushforward of 
$\Sigma$\ by $\Phi_1$. 

Comparing the zero sets of $\Sigma$\ and $\tilde\Sigma$, we see that one of the zeros of $\Sigma$\ has moved from $X$\ to its complement in $\Omega$. 
After repeating this procedure a finite number of times, we are finished.  \cqd

\begin{cor}
\label{tricos} If, in addition, $X$\ is orientable and has dimension two, then the cotangent bundle of the interior, $\tix$, is homeomorphic to 
$\ix\times\re^2$.
\end{cor}
\pf Starting with the section given by Proposition~\ref{tritan}, one may use orientability to get a smooth frame for the cotangent bundle. \cqd

Let $m^\prime:C_0(\ix)\to\cosx$\ denote the pullback of functions under the bundle projection $\six\to\ix$. 
Proposition~\ref{tritan}\ makes it possible to choose a continuous section $\Sigma$\ of the cosphere bundle.
Let then $s:\cosx\to C_0(\ix)$\ denote the \cst-algebra homomorphism $f\mapsto f\circ\Sigma$. It is clear that $s\circ m^\prime$\ is the identity on 
$C_0(\ix)$.

Let $\bix$\ denote the bundle of closed unit balls of the cotangent bundle of the interior. The kernel of the restriction mapping $R:\cobx\to\cosx$\
is homeomorphic to $\cotx$. This observation defines the exact sequence
\begin{equation}
\label{cobundle}
0\to\cotx\to\cobx\to\cosx\to 0.
\end{equation}

\begin{pro}
\label{kths}
If $X$\ is connected and $\bx$\ is not empty, then,
for each $i=0,1$, $K_i(\cosx)$\ is isomorphic to $K_i(C_0(\ix))\oplus K_{1-i}(\cotx)$. Under this isomorphism, the homomorphism induced by $m^\prime$\
corresponds to the canonical injection of $K_i(C_0(\ix))$\ into $K_i(C_0(\ix))\oplus K_{1-i}(\cotx)$, and the connecting mapping $K_i(\cosx)\to K_{1-i}(\cotx)$\ in 
the six-term exact sequence one gets from (\ref{cobundle}) corresponds to the canonical projection of $K_i(C_0(\ix))\oplus K_{1-i}(\cotx)$\ onto 
$K_{1-i}(\cotx)$.
\end{pro}

\pf
Since closed balls can be continuously deformed to their centers, the pullback of functions under the bundle projection $m^{\prime\prime}:C_0(\ix)\to\cobx$\ 
is a homotopy equivalence. Hence, the induced homomorphisms, $m^{\prime\prime}_{*}:K_i(C_0(\ix))\to
K_i(\cobx)$, $i=0,1$, are isomorphisms (\cite{WO}, 6.4.3 and 7.1.6). It is obvious that $R\circ m^{\prime\prime}=m^\prime$. Hence, if we use
the isomorphisms $m^{\prime\prime}_{*}$\ to identify $K_i(\cobx)$\ and $K_i(C_0(\ix))$, the six-term exact sequence associated to (\ref{cobundle}) becomes
\begin{equation}
\label{cclc}
\begin{array}{ccccc}
K_0(\cotx)&\longrightarrow&K_0(C_0(\ix))&{\mathop{\longrightarrow}\limits^{m_{*}^{\prime}}}&K_0(\cosx)\\\big\uparrow&&&&\big\downarrow
\\K_1(\cosx)&{\mathop{\longleftarrow}\limits^{m_{*}^{\prime}}}&K_1(C_0(\ix))&\longleftarrow& K_1(\cotx)
\end{array}.
\end{equation}

We have seen, right after Corollary~\ref{tricos}, that $s\circ m^\prime$\ is the identity. That implies that $s_*$\ is a left inverse for $m^{\prime}_{*}$. 
The cyclic sequence (\ref{cclc}) then becomes two split exact sequences
\[
0\longrightarrow K_i(C_0(\ix)){\mathop{\longleftrightarrow}\limits_{s_{*}}^{m^{\prime}_{*}}}K_i(\cosx)\longrightarrow K_{1-i}(\cotx)\longrightarrow 0,
\]
$i=0,1$. That the above sequences also split on the right and the connecting mappings correspond to projections follows now from algebraic 
generalities (\cite{WO}, 3.1.4, for example). \cqd

\begin{thm}
\label{kthr}
If $X$\ is connected and $\bx$\ is not empty, then, for each $i=0,1$, there is an exact sequence
\begin{equation}
\label{losung}
0\longrightarrow \ker r^i_{*}\oplus K_{1-i}(\cotx)\longrightarrow K_i(\ac/\kc)\longrightarrow \text{{\em Im}}\,r^i_*\longrightarrow 0,
\end{equation}
where $r^i_*:K_{i}(C(X))\to K_i(C(\bx))$\ is the homomorphism induced by the 
restriction to the boundary $r:C(X)\to C(\bx)$.
\end{thm}

\pf
Let us consider the commutative diagram
\begin{equation}
\def\mapup#1{\Big\uparrow\rlap{$\vcenter{\hbox{$\scriptstyle#1$}}$}}
\label{cd}
\begin{array}{ccccccc}
0\longrightarrow&\ic/\kc&\longrightarrow&\ac/\kc&{\mathop{\longrightarrow}\limits^{\pi}}&\ac/\ic&\longrightarrow 0
\\                &\mapup{m}&&\mapup{m}&&\mapup{b}&
\\0\longrightarrow&C_0(\ix)&\longrightarrow&C(X)&{\mathop{\longrightarrow}\limits^{r}}&C(\bx)&\longrightarrow 0
\end{array},
\end{equation}
where $m$\ is the isometric $^*$-homomorphism that maps $f\in C(X)$\ to the class of $\left(\begin{array}{cc}f&0\\0&g\end{array}\right)$, $g$\ denoting the 
restriction to $\bx$\ of $f$; recall that $b$\ was defined in (\ref{bdmp}). Here we do not distinguish between the isomorphic \cst-algebras $\ac/\ic$\ and 
$\img$\ (Corollary~\ref{ami}). 

Let us denote by $\delta$\ and $\exp$\ the index and exponential mappings associated to the top exact sequence in (\ref{cd}), and 
by $\delta^0$\ and $\exp^0$\ the index and exponential mappings associated to the bottom one. By the naturality of the connecting mappings
\cite{RLL}, we get from (\ref{cd}) the two commutative diagrams
\begin{equation}
\label{cdt}
\def\mapup#1{\Big\uparrow\rlap{$\vcenter{\hbox{$\scriptstyle#1$}}$}}
\begin{array}{ccccccc}
K_1(\ac/\ic)&{\mathop{\longrightarrow}\limits^{\delta}}&K_0(\ic/\kc)&&K_0(\ac/\ic)&{\mathop{\longrightarrow}\limits^{\text{exp}}}&K_1(\ic/\kc)\\
\mapup{b_{*}}&&\mapup{m_{*}}&\ \text{and}\ &\mapup{b_{*}}&&\mapup{m_{*}}\\
K_1(C(\bx))&{\mathop{\longrightarrow}\limits^{\delta^{0}}}&K_0(C_0(\ix))&&K_0(C(\bx))&{\mathop{\longrightarrow}\limits^{\text{exp}^{0}}}&K_1(C_0(\ix))
\end{array}.
\end{equation}

The \cst-algebra homomorphism $m^\prime$, defined after Corollary~\ref{tricos}, is the composition of the isomorphism $j:\ic/\kc\to\cosx$\ of 
Theorem~\ref{tker}\ with 
$m$, defined in (\ref{cd}). Inserting the isomorphisms $j_*:K_i(\ic/\kc)\to K_i(\cosx)$\ into the upper right corners of the diagrams in (\ref{cdt}), and 
denoting by $\delta^\prime$\ and $\exp^\prime$\ the compositions of $j_*$\ with $\delta$\ and $\exp$, we get:
\begin{equation}
\label{cdtl}
\def\mapup#1{\Big\uparrow\rlap{$\vcenter{\hbox{$\scriptstyle#1$}}$}}
\begin{array}{ccccccc}
K_1(\ac/\ic)&{\mathop{\longrightarrow}\limits^{\delta^{\prime}}}&K_0(\cosx)&&K_0(\ac/\ic)&{\mathop{\longrightarrow}\limits^{\text{exp}^{\prime}}}&K_1(\cosx)\\
\mapup{b_{*}}&&\mapup{m^\prime_{*}}&\text{and}&\mapup{b_{*}}&&\mapup{m^\prime_{*}}\\
K_1(C(\bx))&{\mathop{\longrightarrow}\limits^{\delta^{0}}}&K_0(C_0(\ix))&&K_0(C(\bx))&{\mathop{\longrightarrow}\limits^{\text{exp}^{0}}}&K_1(C_0(\ix))
\end{array}.
\end{equation}

In other words, modulo isomorphisms and split injections, the connecting mappings in the six-terms exact sequences associated to the two horizontal short exact sequences in (\ref{cd}) are the same. We show, next, that this reduces the computation of the connecting mappings in the six-term exact sequence associated (with the use of the isomorphisms $b_*$\ and $j_*$) to $0\to\ic/\kc\to\ac/\kc\to\ac/\ic\to 0$, 
\begin{equation}
\label{ciclic}
\def\mapup#1{\Big\uparrow\rlap{$\vcenter{\hbox{$\scriptstyle#1$}}$}}
\def\mapdown#1{\Big\downarrow\rlap{$\vcenter{\hbox{$\scriptstyle#1$}}$}}
\begin{array}{ccccc}
K_0(C_0(\six))&\longrightarrow&K_0(\ac/\kc)&\longrightarrow&K_0(C(\bx))\\
\mapup{\delta^{\prime\prime}}&&&&\mapdown{\text{exp}^{\prime\prime}}\\
K_1(C(\bx))&\longleftarrow&K_1(\ac/\kc)&\longleftarrow&K_1(C_0(\six))
\end{array},
\end{equation}
to a purely topological operation. In the diagram, above, we have denoted $\delta^\prime\circ b_*$\ and $\text{exp}^\prime\circ b_*$\ by $\delta^{\prime\prime}$\
and $\exp^{\prime\prime}$, respectively. 

Taking quotients by the kernels of the upper-left and lower-right horizontal arrows, and restricting the ranges of the other horizontal arrows in (\ref{ciclic}), we 
obtain the exact sequences
\begin{equation}
\label{detalhe1}
0\longrightarrow\frac{K_0(C_0(\six))}{\text{Im}\,\delta^{\prime\prime}}\longrightarrow K_0(\ac/\kc)\longrightarrow\ker{\text{exp}^{\prime\prime}}\longrightarrow 0,
\end{equation}
and
\begin{equation}
\label{detalhe2}
0\longrightarrow\frac{K_1(C_0(\six))}{\text{Im}\,\text{exp}^{\prime\prime}}\longrightarrow K_1(\ac/\kc)\longrightarrow\ker{\delta^{\prime\prime}}\longrightarrow 0.
\end{equation}

A similar argument with quotients, applied to the cyclic sequence associated to the bottom exact sequence in (\ref{cd}),
\begin{equation}
\def\mapup#1{\Big\uparrow\rlap{$\vcenter{\hbox{$\scriptstyle#1$}}$}}
\def\mapdown#1{\Big\downarrow\rlap{$\vcenter{\hbox{$\scriptstyle#1$}}$}}
\label{cyclic}
\begin{array}{ccccc}
K_0(C_0(\ix))&\longrightarrow&K_0(C(X))&{\mathop{\longrightarrow}\limits^{r_*}}&K_0(C(\bx))\\
\mapup{\delta^0} &&&&\mapdown{\exp^0}\\
K_1(C(\bx))&{\mathop{\longleftarrow}\limits^{r_*}}&K_1(C(X))&\longleftarrow&K_1(C_0(\ix))
\end{array},
\end{equation}
gives the isomorphisms 
\begin{equation}
\label{40}
\frac{K_0(C_0(\ix))}{\text{Im}\,\delta^0}\simeq\ker r_*^0\ \ \ \text{and}\ \ \ \frac{K_1(C_0(\ix))}{\text{Im}\,\text{exp}^0}\simeq\ker r_*^1.
\end{equation}

By (\ref{cdtl}), $\text{Im}\,\delta^{\prime\prime}=\text{Im}(m^\prime_*\circ\delta^0)$\ and $\text{Im}\,\text{exp}^{\prime\prime}=\text{Im}(m_*^\prime\circ\text{exp}^0)$. 
Proposition~\ref{kths}\ then implies
\begin{equation}
\label{41}
\frac{K_0(C_0(\six))}{\text{Im}\,\delta^{\prime\prime}} \simeq K_1(\cotx)\oplus \frac{K_0(C_0(\ix))}{\text{Im}\,\delta^0},
\end{equation}
and
\begin{equation}
\label{42} 
\frac{K_1(\cosx)}{\text{Im}\,\text{exp}^{\prime\prime}} \simeq K_0(\cotx)\oplus \frac{K_1(C_0(\ix))}{\text{Im}\,\text{exp}^0}.
\end{equation}
Hence, the groups at the left in (\ref{losung}) are isomorphic to those in (\ref{detalhe1}) and (\ref{detalhe2}) by (\ref{40}), (\ref{41}) and (\ref{42}). 

Also by (\ref{cdtl}), $\ker\delta^{\prime\prime}=\ker\delta^0$, and $\ker\text{exp}^{\prime\prime}=\ker\text{exp}^0$. Hence, the groups at the right of (\ref{losung})
are isomorphic to those of (\ref{detalhe1}) and (\ref{detalhe2}) by the exactness of (\ref{cyclic}). \cqd.

It is well known that the K-groups of $C_0(Y)$\ are finitely generated, for any manifold $Y$. That follows by induction, using triangularizability,
starting from the fact that $K_0(C_0(Y))\simeq\co$\ and $K_1(C_0(Y))=0$\ when $Y$\ is a point. 

\begin{cor}
\label{free}
If, in addition to the hypothesis of Theorem~\ref{kthr}, the K-groups of $C(\bx)$\ have no torsion, then 
$K_i(\ac/\kc)\simeq K_i(C(X))\oplus K_{1-i}(\cotx)$, $i=0,1$.
\end{cor}
\pf By our previous remark, each $K_i(C(\bx)$, $i=0,1$, is finitely generated. The hypothesis then implies that they are free, and so are
their subgroups $\text{Im}\,r_*^i$.

If $0 \to A \to B \to C \to 0$ is a short exact sequence of abelian groups with $C$ free, then $B$ is
isomorphic to $A \oplus C$. Our claim follows from this fact applied to (\ref{losung}) and to  $0 \to \ker r^i_* \to C(X) \to \text{Im}\,r ^i_* \to 0$.\cqd
 
In Section \ref{superficie}, we apply Corollary~\ref{free}\ to orientable surfaces. We end this section showing how one can get the K-theory of $\ac$\ and 
$\as$\ from the K-theory of $\ac/\kc$. The result about $K_1$\ in Proposition~\ref{kac}\ is improved in Corollary~\ref{elmar}.

\begin{pro}\label{kac}
The projection $\ac\to\ac/\kc$ induces isomorphisms from $K_0(\ac)$\ to $K_0(\ac/\kc)$, and from $K_1(\ac)$\ to the kernel of the Fredholm-index 
mapping $K_1(\ac/\kc)\to\z$, which is surjective.
\end{pro}
\pf As in the proof of Lemma~\ref{ktht}, it suffices to prove that there exist an integer $k$\ and a $k$-by-$k$\ matrix with entries in $\ac$\ 
which, regarded as an \op\ on $\h^k$, is Fredholm and has index one.

It follows from Fedosov's index formula, as in the proof of Theorem~5.18 in \cite{LM}, that there exists a $k$-by-$k$\ matrix $S$\ of zero-order \psd\ \ops\ 
on $\bx$\ defining an index-one Fredholm \op\ on $(H^{-\frac{1}{2}}(\bx))^k$. 
The Green \op\ $\left(\begin{array}{cc}I&0\\0&S\end{array}\right)$, $I$\ denoting the identity on $C^\infty(X;\co^k)$, is then a Fredholm \op\ of index
one on $\h^k$. \cqd

\begin{pro}\label{kas}
The injection of $\as$\ into $\ac$\ induces isomorphisms between $K_i(\as)$\ and $K_i(\ac)$, $i=0,1$.
\end{pro}
\pf 
By \cite{S2}, Corollary 4.11, $\as$\ can be given the structure of a Fr\'echet $^{*}$-algebra, such that the embedding of $\as$\ in $\ac$\ is continuous.
Moreover, $\as$\ contains the inverses of all its elements which are invertible in $\ac$. In particular, the set of invertibles in $\as$\ is open, with 
respect to that Fr\'echet topology. Then, by \cite{W}, p.\ 115, the inversion
is continuous. Hence, the Cauchy integrals that give the holomorphic functional 
calculus converge also in $\as$. Being closed under the holomorphic functional calculus, $\as$\ has the same K-theory as $\ac$\ 
(\cite{Bost}, Th\'eor\`eme A.2.1). 
\cqd

%
%
%
%
%
%
\section{A better result for $K_1$}\label{index}

Throughout this section, we assume that $X$\ is connected and that $\bx$\ is not empty, in order to be able to apply Proposition~\ref{kths}\ and the proof 
of Theorem~\ref{kthr}.

As a first step, let us show how the principal symbol is related to the index mapping in the six-term exact sequence associated to (\ref{cobundle}). 
We need topological K-theory and refer to the first section of \cite{Bl}\ for definitions and notation. 

Let $Y$\ be a locally compact Hausdorff space, $Z$\ a closed subspace of $Y$, and $U=Y\setminus Z$. Let $\imath:K_0(C_0(U))\to K(Y,Z)$\ denote the composition of
the isomorphism from $K(U)$\ to $K(Y,Z)$\ given in \cite{Bl}, 1.5.1, with the canonical isomorphism from $K_0(C_0(U))$\ to $K(U)$. 

We learned the following proof from Ruy Exel.

\begin{pro}
\label{ruy}
Let $\delta:K_1(C_0(Z))\to K_0(C_0(U))$\ be the index mapping in the six-term exact sequence associated to $0\to C_0(U)\to C_0(Y)\to C_0(Z)\to 0$. Then
$\imath\circ\delta:K_1(C_0(Z))\to K(Y,Z)$\ maps the class of an invertible $k$-by-$k$\ matrix $u\in M_k(C_0(Z)^+)$, $C_0(Z)^+=C(Z^+)=C(Z\cup\{+\})$, 
to the class of the triple $(E^k,E^k,u)$, where $E^k$\ denotes the rank-$k$\ trivial bundle and $u$\ is regarded as an isomorphism on the restriction of 
$E^k$\ to  $Z$\ in the obvious way.
\end{pro}
\pf
We may suppose, without loss of generality, that $u(+)$\ is the identity.
Let $w\in M_{2k}(C_0(Y)^+)$\ be an invertible coinciding with $u\oplus u^{-1}\in M_{2k}(C_0(Z)^+)$\ on $Z$\ and such that $w(+)$\ is the identity. By 
definition (\cite{WO}, 8.1.1), 
$\delta([u])=[wp_kw^{-1}]-[p_k]$, where $p_k$\ is the $2k$-by-$2k$\ matrix with $1$\ on the first $k$\ entries of the diagonal and zero elsewhere. The above 
mentioned canonical isomorphism maps an element $[p]-[q]\in K_0(C_0(U))$, $p$\ and $q$\ idempotents in $M_l(C_0(U)^+)$, to the element
of $K(U)\simeq K(U^+,+)$\ defined by the triple $(\text{Im}\,p,\text{Im}\,q,\alpha)$, with $\text{Im}\,p$\ denoting the vector bundle 
$\{(x,v)\in U^+\times\co^l;v\in\text{Im}\,p(x)\}$, and $\alpha$\ being any isomorphism between $\text{Im}\,p(+)$\ and $\text{Im}\,q(+)$\ 
({\em any}\ homomorphism because the equivalence class of this triple depends only on the homotopy class of $\alpha$, by \cite{K}, II.2.15, and the 
complex linear group $\text{Gl}_k(\co)$\ is connected; this also justifies the first statement in this proof). 
Since $wp_kw^{-1}$\ is equal to $p_k$\ at the infinite point $+$, the fibers
of $\text{Im}\,wp_kw^{-1}$\ and $\text{Im}\,p_k$\ are canonically isomorphic there. Hence,  viewed as an element of $K(U)$, $\delta([u])$\ is the class of 
the triple $(\text{Im}\,wp_kw^{-1},\text{Im}\,p_k,\text{id})$, where $\text{id}$\ denotes the identity mapping.

Let $\phi:Y^+\to U^+$\ be the identity on $U$\ and map all other points to $+$. The isomorphism from $K(U)$\ to $K(Y,Z)$\ defined in \cite{Bl}, 1.5.1, is the 
homomorphism contravariantly induced by $\phi$, viewed as a morphism between the compact pairs $(Y^+,Z^+)$\ and $(U^+,+)$. Hence, it maps $\delta([u])$\ 
to the class of the triple $(\text{Im}\,wp_kw^{-1},\text{Im}\,p_k,\text{id})$\ (all that $\phi^*$\ does is to regard the projections defining the triple
as functions on $Y^+$\ which are constant over $Z^+$). 

The mapping $f:\text{Im}\,p_k\to\text{Im}\,wp_kw^{-1}$, $f(x,v)=(x,w(x)v)$, is a vector bundle isomorphism. Moreover, $f$\ and the identity mapping
on $\text{Im}\,p_k$\ intertwine $u\oplus u^{-1}$\ and the identity. The triple $(\text{Im}\,wp_kw^{-1},\text{Im}\,p_k,\text{id})$\ is therefore equivalent
to $(\text{Im}\,p_k,\text{Im}\,p_k,u\oplus u^{-1})$, which is obviously equivalent to $(E^k,E^k,u)$. \cqd

We have seen in Proposition~\ref{kths}\ that $K_1(\cosx)$\ is isomorphic to $K_1(C_0(\ix))\oplus K_0(\cotx)$, with the canonical projection onto 
$K_0(\cotx)$\ corresponding to the index mapping $\delta^1:K_1(\cosx)\to K_0(\cotx)$\ associated to (\ref{cobundle}). Applying Proposition~\ref{ruy}, with 
$Y=\bix$\ and $Z=\six$, we may give another description of $\delta^1$, now regarded as a mapping from $K_1(\ic/\kc)$\ to $K_0(\cotx))$\ 
(by Theorem~\ref{tker}). Let $[[A]]\in K_1(\ic/\kc)$\ be given, $[[A]]$\ denoting the $K_1$-class of the class $[A]\in(\ic/\kc)^+\subset\ac/\kc$\ of a 
Fredholm \op\ $A\in M_k(\ic\oplus\co)$\ with principal symbol $\sig(A)$\ (we denote by ${\mathfrak C}^+$\ the unitization of a \cst-algebra ${\mathfrak C}$).
$\delta^1([[A]])$\ is then equal to the class of the triple $(E^k,E^k,\sig(A))$\ in $K(\bix,\six)\simeq K(\tix)\simeq K_0(\cotx)$.

In our language, part of the content of \bdm's index theorem is that $\delta^1:K_1(\cosx)\to K_0(\cotx)$\ factors through 
$K_1(\ac/\kc)$. More precisely, let us denote by $\ind:K_1(\ac/\kc)\to K_0(\cotx)$\ the composition of the canonical isomorphism between $K(\tix)$\ and
$K_0(\cotx)$\ with the homomorphism defined in \cite{B2}, Theorem 5.21 (see also \cite{RS}, 3.2.2.4, Theorem~1). We then have:

\begin{lem}\label{comuta}
With $i_*$\ denoting the homomorphism induced by the inclusion of $\ic/\kc$\ into $\ac/\kc$, one has $\ind\circ i_*=\delta^1$.
\end{lem}
\pf
Let $x=[[A]]\in K_1(\ic/\kc)$\ be given, $A$\ as above. Since the scalar part of $A$\ is invertible and $\text{Gl}_k(\co)$\ is connected, we may suppose that
$A=B+I$, for some $B\in M_k(\ic)$\ and $I$\ the $k$-by-$k$\ identity matrix. By definition of $\ic$, there is a sequence $B_j\to B$, 
$B_j=\left(\begin{array}{cc}\varphi P\psi+G&K\\T&S\end{array}\right)\in\is$, with $P$, $G$, etc, denoting $k$-by-$k$\ matrices of \ops\ as those in 
(\ref{dis}). For sufficiently large $j$, then, $x=[[I+B_j]]$\ (by \cite{RLL}, 2.1.11). Since $\left(\begin{array}{cc} G&K\\T&S\end{array}\right)$\ is compact,
$x=\left[\left[\left(\begin{array}{cc} \varphi P\psi+I&0\\0&I\end{array}\right)\right]\right]$, with $I$\ denoting the identity operators on 
$\cix$\ and on $\cibx$. This \op\ is in the form \cite{B2}, (5.21)-(3). According to \bdm's prescription, $\ind(x)$\ corresponds to the class in 
$K(\tix)\simeq K(\bix^+,\six^+)$\ determined by the principal symbol of $\varphi P\psi+I$\ (the contribution from the identity on $\cibx$\ vanishes).
This proves the lemma, by our comments after the proof of Proposition~\ref{ruy}. \cqd

Using Lemma~\ref{comuta}, Proposition~\ref{kths}, and the proof of Theorem~\ref{kthr}, we then obtain the commutative diagram:

\begin{equation}
\def\mapup#1{\Big\uparrow\rlap{$\vcenter{\hbox{$\scriptstyle#1$}}$}}
\def\mapnw#1{\nwarrow\rlap{$\vcenter{\hbox{$\scriptstyle#1$}}$}}
\label{cdbdm}
\begin{array}{cccccc}
K_0(\cotx)&&&&&\\
\mapup{\delta^1}&\mapnw{\ind}&&&&\\
K_1(C_0(\ix))\oplus K_0(\cotx)&{\mathop{\longrightarrow}\limits^{i_*}}&K_1(\ac/\kc)&{\mathop{\longrightarrow}\limits^{\pi_*}}&K_1(\ac/\ic)
&{\mathop{\longrightarrow}\limits^{\delta^\prime}}\\
\mapup{m^\prime_*}&&\mapup{m_*}&&\mapup{b_*}&\\
K_1(C_0(\ix))&{\mathop{\longrightarrow}\limits^{i^\prime_*}}&K_1(C(X))&{\mathop{\longrightarrow}\limits^{r_*}}&K_1(C(\bx))
&{\mathop{\longrightarrow}\limits^{\delta^0}}
\end{array}.
\end{equation}

\begin{lem}
\label{aladim}
In the diagram, above, we have $\ind\circ m_*=0$.
\end{lem}
\pf We want to show that $\ind(x)=0$, whenever $x\in K_1(\ac/\kc)$\ is of the form 
$x=[[A]]$, with $A=\left(\begin{array}{cc} f&0\\0&f|_{\bx}\end{array}\right)$, for some invertible $f\in M_k(C(X))$. 
Such an $A$\ is in the form of \cite{RS}, 3.2.2.4, Theorem 1, (iii). Indeed, the symbol $f$\ of its upper left corner is independent of the covariable not
only on a neighborhood of the boundary, but over all $X$. The class of the triple $(E^k,E^k,f)$\ in $K(B^*X,S^*X\cup T^*X|_{\bx})\simeq K(\tix)$\ is zero, 
because $f$\ defines a bundle isomorphism of $E^k=B^*X\times\co^k$\ onto itself. For the same reason, the element of $K(T^*\bx)$\ determined by $f|_{\bx}$\ 
also vanishes. Hence, $\ind(x)=0$. \cqd

\begin{thm}
\label{nosso} If $X$\ is connected and $\bx$\ is not empty, then $K_1(\ac/\kc)$\ is isomorphic to $K_1(C(X))\oplus K_0(\cotx)$. More precisely, in the 
diagram (\ref{cdbdm}), $m_*$\ and the restriction of $i_*$\ to $K_0(\cotx)$\ are injective, and $K_1(\ac/\kc)=m_*(K_1(C(X))\oplus i_*(K_0(\cotx))$.
\end{thm}
\pf
Given $x\in K_1(\ac/\kc)$, we have $\delta^0(b_*^{-1}(\pi_*(x)))=0$, because $\delta^\prime(\pi_*(x))=0$, $m_*^\prime$\ is injective on $K_0(C_0(\ix))$, and 
the left diagram in (\ref{cdtl}) commutes. It is therefore possible to choose $y\in K_1(C(X))$\ such that $r_*(y)=b_*^{-1}(\pi_*(x))$. Since 
$\pi_*(x-m_*(y))=0$, there exists $z_1\oplus z_2\in K_1(C_0(\ix))\oplus K_0(\cotx)$\ such that $i_*(z_1\oplus z_2)=x-m_*(y)$.

In Proposition \ref{kths}, we showed that $m_*^\prime$\ is the canonical injection of $K_1(C_0(\ix))$\ into $K_1(C_0(\ix))\oplus K_0(\cotx)$. The 
commutativity of the lower left subdiagram in (\ref{cdbdm}) then implies that $x=m_*(i_*^\prime(z_1))+m_*(y)+i_*(z_2)$. This proves that 
$m_*(K_1(C(X))+i_*(K_0(\cotx))=K_1(\ac/\kc)$. 

To show that the intersection of $m_*(K_1(C(X))$\ and $i_*(K_0(\cotx))$\ is $0$, let us suppose that $m_*(y)=i_*(0\oplus z)$, for some $y\in K_1(C(X))$\ and 
$z\in K_0(\cotx)$. Since $\ind(m_*(y))=0$\ (Lemma~\ref{aladim}), $z=\ind(i_*(0\oplus z))=\ind(m_*(y))=0$. 

The existence of the homomorphism $\ind$\ implies at once that $i_*$\ restricted to $K_0(\cotx)$\ is injective: $\ind(i_*(0\oplus z))=z$.

To prove that $m_*$\ is injective, let an invertible $f\in M_k(C(X))$\ be given, such that $A=\left(\begin{array}{cc} f&0\\0&f|_{\bx}\end{array}\right)$\ 
can be connected to the identity by a continuous path of Fredholm \ops\ in $M_k(\ac)$. Using that $[A]\mapsto\sig(A)$\ defines a continuous mapping from 
$M_k(\ac/\kc)$\ to $M_k(C(S^*X))$\ (this follows from (\ref{rsc})), we then get a homotopy of invertibles in $M_k(C(S^*X))$, between the $k$-by-$k$\ identity
matrix $I$\ and $f$\ (regarded as a function on $S^*X$\ independent of the covariable). In Proposition~\ref{tritan}, we showed that there exists a continuous 
section $\Sigma$\ of $S^*X$. By composition with $\Sigma$, any homotopy in $M_k(C(S^*X))$\ defines a homotopy in $M_k(C(X))$. 
Hence, there is a homotopy of invertibles in $M_k(C(X))$\ connecting $f$\ and $I$. The class defined by $f$\ in $K_1(C(X))$\ therefore vanishes.  \cqd

The following two corollaries follow immediately from Lemma~\ref{aladim}\ and Theorem~\ref{nosso}. 

\begin{cor}\label{it}
With respect to the isomorphism of Theorem~\ref{nosso}, $\ind$\ corresponds to the canonical projection from $K_1(C(X))\oplus K_0(\cotx)$\ onto
$K_0(\cotx)$.
\end{cor}

\begin{cor}\label{savin}
$0 \longrightarrow K_1(C(X)) {\mathop{\longrightarrow}\limits^{m_{*}}} K_1(\ac/\kc) {\mathop{\longrightarrow}\limits^{\ind}} K_0(\cotx)\longrightarrow 0$\ 
is exact.
\end{cor}

This section was inspired by conversations with Anton Savin about \bdm's index theorem, at conferences in Potsdam and B\c edlewo. Corollary~\ref{savin}\ is 
his conjecture. 

Next we show that also $K_1(\ac)$\ is topologically determined. We are going to use that there exists a mapping $\chi:K_0(\cotx)\to\z$\ (the 
{\em topological index}) such that $\chi\circ\ind$\ gives the Fredholm index (\cite{B2}, Section 5.8; \cite{RS}, 3.2.2.3, 3.2.2.4). 

\begin{cor}\label{elmar}
$K_1(\ac)$\ and $K_1(\as)$\ are isomorphic to $K_1(C(X))\oplus\ker\chi$.
\end{cor}
\pf
It follows from \bdm's index theorem, quoted above, and from our Proposition~\ref{kac}, that $K_1(\ac)$\ is isomorphic to $\ker(\chi\circ\ind)$. 
Theorem~\ref{nosso}\ and Corollary~\ref{it}\ imply that an arbitrary element of $K_1(\ac/\kc)$\ is of the form $m_*(x)\oplus i_*(y)$, $x\in K_1(C(X))$\ and
$y\in K_0(\cotx)$, and that $\ind(m_*(x)\oplus i_*(y))=y$. Then it is obvious that $\chi\circ\ind(m_*(x)\oplus i_*(y))=0$\ \ifoi\ $y\in\ker\chi$. 
In Theorem~\ref{nosso}, we also proved that $m_*$\ is injective. \cqd

%
%
%
%
%
%
\section{Orientable Surfaces}\label{superficie}

Throughout this section, $X$\ denotes a connected orientable two-dimensional manifold with nonempty boundary $\bx$. The genus of $X$\ is denoted by $g$\ and 
the number of connected components of $\bx$\ by $m$. 

It is probably well known (we thank Thomas Schick for this information) that any manifold like our $X$\ can be continuously deformed to the 
union of $q$\ circles with one point in common, $q=2g+m-1$. That already implies that $K_0(C(X))\simeq\z$\ and $K_1(C(X))\simeq\z^{q}$. 
In the proof of Proposition~\ref{ktx}, we give a precise description of that deformation, which is also used to prove Proposition~\ref{ktix}. 

\begin{pro}
\label{ktx} There exist $V\in X$, closed curves $C_1,\cdots,C_q$\ in $X$, $q=2g+m-1$, such that $C_j\cap C_k=\{V\}$\ if $j\neq k$, and an isomorphism from 
$\z^{q}$\ to $K_1(C(X))$, which maps each element $e_k$\ of the canonical basis to the class of a unitary in $C(X)$\ equal to one on $C_j$, if 
$j\neq k$, and with winding number one on $C_k$. Moreover, $K_0(C(X))=[1]\cdot\z$, where $[1]$\ denotes the class of the function 
identical to $1$\ on $X$.
\end{pro}

\pf
Any closed orientable surface of genus $g$\ is homeomorphic to a polygon of $4g$\ sides, identified in pairs, all vertices corresponding to the same point in
the manifold (\cite{F}, Section 17b). Orientability implies that, if $l$\ and $l^\prime$\ are two identified sides, then the polygon is to the left of $l$\ if 
and only if it is to the right of $l^\prime$\ (assuming, of course, that the parametrizations of $l$\ and $l^\prime$\ are the same, with respect to the 
identification). Our surface X is then homeomorphic to such a polygon with sides identified, with $m$\ disjoint open disks, $D_1,\cdots,D_m$, 
removed from its interior. Let us choose a side $l=[P,Q]$\ in this polygon and draw $m-1$\ curves, $c_1,\cdots,c_{m-1}$, all going from $P$\ to $Q$, so that 
$D_1$\ is between $l$\ and $c_1$; $D_k$\ is between $c_{k-1}$\ and $c_k$, $k=2,\cdots,m-1$; and $D_m$\ is between $c_{m-1}$\ and the remaining sides. 

Each curve $c_k$\ and each pair of identified sides in the polygon correspond to circles (closed curves) in $X$\ that do not intersect $\bx$. Let us denote
by $C_k$\ the circles obtained from $c_k$, $k=1,\cdots,m-1$; by $C_m$\ the circle that comes from $l$\ and its pair; and by $C_{m+1},\cdots,C_{q}$\ the circles
that come from the other sides of the polygon. Any two 
among these $q$\ circles meet in exactly one point, $V$, the equivalence class of the vertices of the polygon. Let us denote by $Y$\ the union of the 
circles $C_1,\cdots,C_{q}$. Gradually enlarging the disks $D_1,\cdots,D_m$, without crossing any of the $C_k$'s, but with their boundaries eventually adhering
to them, one proves that $X$\ is homotopically equivalent to $Y$. 

Looking at what this deformation does to continuous functions on $X$, one proves that the restriction mapping $R:C(X)\to C(Y)$\ is a homotopy equivalence. 
Hence, we get the isomorphisms
\begin{equation}
\label{rstar}
R_*:K_i(C(X))\to K_i(C(Y)),\ i=0,1. 
\end{equation}

For any circle $C$\ containing a point $V$, let us denote by $C_0(C\backslash V)$\ the algebra of continuous functions on $C$\ that vanish at $V$.
We are going to use that $K_0(C_0(C\backslash V))=0$; and that $K_1(C_0(C\backslash V))\simeq\z$, with isomorphism given by the winding number. 
Moreover, there exists a generator of $K_1(C_0(C\backslash V))$\ which is equal to one at $V$. 

Let us consider the exact sequence
\begin{equation}
\label{stun}
0 \longrightarrow\bigoplus_{k=1}^{q}C_0(C_k\backslash V)
\smash{\mathop{\longrightarrow}\limits^i} C(Y)\smash{\mathop{\longrightarrow}\limits^\pi}\co\longrightarrow 0,
\end{equation}
where $i$\ denotes the inclusion mapping, and $\pi$\ evaluation at $V$. Using that $K_0(\co)\simeq\z$\ and $K_1(\co)=0$, we get from (\ref{stun}):
\begin{equation}
\label{kty}
\begin{array}{ccccl}
0&\longrightarrow&K_0(C(Y))&\longrightarrow&\z\\
\big\uparrow&&&&\big\downarrow\\
0&\longleftarrow&K_1(C(Y))&\longleftarrow&\z^{q}
\end{array}.
\end{equation}

The exponential mapping in (\ref{kty}) is the zero map, since the upper-right horizontal arrow  maps $[1]$\ to $1\in\z$. 
We then get that $K_0(C(Y))$\ is isomorphic to $\z$, and the statement about $K_0$\ follows from (\ref{rstar}). We also get that the inclusion mapping in 
(\ref{stun}) induces an isomorphism $i_*$\ from $\oplus_kK_1(C_0(C_k\backslash V))\simeq\z^{q}$\ to $K_1(C(Y))$. The composition 
$R_*^{-1}\circ i_*$\ gives 
the other isomorphism for which we were looking. \cqd

$K_0(C(\bx))$\ and $K_1(C(\bx))$\ are both isomorphic to $\z^m$. The generators of $K_0(C(\bx))$\ are the classes of the functions equal to one on the $k$-th 
connected component of $\bx$\ (the boundary of $D_k$, as defined above), and zero on the others. We choose the isomorphism from $K_1(C(\bx))$\ to $\z^m$\ 
given by the winding number, with respect to the orientation of $\bx$\ induced by the orientation of $X$.

Let us now consider the exact sequence
$0\to C_0(\ix)\to C(X) \to C(\bx)\to 0$\ induced by the restriction to the boundary.
Detailed information about all the homomorphisms in the corresponding six-term exact sequence,
\begin{equation}
\label{ciclica}
\begin{array}{ccccc}
K_0(C_0(\ix))&\longrightarrow&\z&\longrightarrow&\z^m\\
\big\uparrow &&&&\big\downarrow\\
\z^m&\longleftarrow&\z^{q}&\longleftarrow&K_1(C_0(\ix))
\end{array},
\end{equation}
is given in the following proposition.

\begin{pro}
\label{ktix}
$K_0(C_0(\ix))$\ is isomorphic to $\z$\ and $K_1(C_0(\ix))$\ is isomorphic to $\z^q$. The index mapping in (\ref{ciclica}) is surjective and, with respect to 
the previously defined isomorphisms, has kernel equal to $\{(j_1,\cdots,j_m);\sum j_k=0\}$. Moreover, the exponential mapping has kernel generated by 
$(1,\cdots,1)$\ and the lower-right arrow in (\ref{ciclica}) has image isomorphic to $\z^{2g}$, generated by $e_{m+1},\cdots,e_q$\ and $(e_1+\cdots+e_m)$.
\end{pro}

\pf
The upper-right horizontal arrow in (\ref{ciclica}) maps $1\in\z$\ to the nonzero element $([1],\cdots,[1])$\ of $K_0(C(\bx))$. 
That means that the upper-left arrow is the zero mapping and that the kernel of the exponential is $(1,\cdots,1)\cdot\z$.

Paying close attention to how the connected components of $\bx$\ touch the $C_k$'s, at the end of the deformation described in the proof of 
Proposition~\ref{ktx}, we see that the lower-left horizontal arrow in (\ref{ciclica}) maps $e_k$\ to $e_{k+1}-e_k$, $k=1,\cdots,m-1$; $e_m$\ to $e_1-e_m$;
and the remaining $e_{m+1},\cdots e_q$\ to zero. We then get that the kernel of the index mapping is what we stated, and that the kernel of the lower-left 
arrow in (\ref{ciclica}) is generated by $e_{m+1},\cdots,e_q$\ and $(e_1+\cdots+e_m)$. Moreover, the quotient of $\z^m$\ by the kernel of the index mapping 
is isomorphic to $\z$, what proves our statement about $K_0(C_0(\ix))$. Finally, because $K_1(C_0(\ix))$\ sits in the 
middle of the exact sequence
\[
0\longrightarrow\frac{\z^m}{(1,\cdots,1)\cdot\z}\longrightarrow K_1(C_0(\ix))\longrightarrow\z^{2g}\longrightarrow 0,
\]
it is free and finitely generated, hence isomorphic to $\z^{2g+m-1}$. \cqd

It now follows from Corollary~\ref{tricos}\ and from Bott periodicity that $K_0(\cotx)\simeq\z$\ and $K_1(\cotx)\simeq\z^q$. From this, Corollary~\ref{free}\
and Proposition~\ref{ktx}, we get:

\begin{cor}
\label{surf}
$K_0(\ac/\kc)$\ and $K_1(\ac/\kc)$\ are both isomorphic to $\z^{2g+m}$.
\end{cor}

Propositions \ref{kac}\ and \ref{kas}\ then imply:

\begin{cor}
\label{abs}
$K_0(\ac)\simeq K_0(\as)\simeq\z^{2g+m}$\ and $K_1(\ac)\simeq K_1(\as)\simeq\z^{2g+m-1}$.
\end{cor}

%
%
%
%
%
%
\section{Composition Sequence}\label{compseq}

In this section we return to our initial assumption: $X$\ is an arbitrary compact manifold with boundary.
We recall that $\gs$\ and $\gc$\ were defined after (\ref{cmpsq}). It is obvious that $\gs=\ker\sigma$. Moreover, $\gs$\ contains all integral \ops\ with
smooth kernel; hence, $\kc\subset\gc$. 

\begin{lem}
\label{novogoh}
There exists a positive constant $c$, determined only by $X$, such that
\begin{equation}
\label{tragog}
\inf_{G^\prime\in\gc}||A+G^\prime||\leq c||\sig(A)||,
\end{equation}
for every $A\in\ac$.
\end{lem}
\pf It is enough to prove (\ref{tragog}) for $A=\pp\in\ac_{11}$, with $P$\ denoting a classical zero-order \psd\ \op\ with the \tp\ on $X$. 
Moreover, since $\kc\subset\gc_{11}$, it suffices to prove that $\inf_{C\in\kc}||\pp+C||\leq c||\sigma(P)||$.

We need to distinguish between the given $P$, regarded as a \psd\ \op\ on $\ix$, and an extension $\tilde P$\ of $P$\ to $\Omega$. 
It follows from the classical norm estimate for \psd\ \ops\ on closed manifolds (quoted between (\ref{lcl}) and (\ref{nev})\,) that, for any $\delta >0$, 
there is a compact $\tilde C$\ on $L^2(\Omega)$\ such that $||\tilde P+\tilde C||<(1+\delta)||\sigma(\tilde{P})||$. If $C$\ is the compact \op\ on $L^2(X)$\ 
obtained from $\tilde C$\ by truncation (i.e., $C$\ is the restriction to $X$, composed with $\tilde C$, composed with zero-extension), then it is obvious 
that $||\pp+C||\leq||\tilde P+\tilde C||$. Hence, the lemma will be proven if we show that the choice of $\tilde P$\ can be made so that
$||\sigma(\tilde P)||\leq c||\sigma(P)||$, for some constant $c$\ depending only on $X$. But that follows from Seeley's extension \cite{Seee}, as in 
\cite{S2}, Lemma~2.3. \cqd

A necessary and sufficient condition for a polyhomogeneous \psd\ \op\ to have the \tp\ was given by \bdm\ \cite{B1}\ (see also \cite{H}, Theorem 18.2.15). For 
zero-order \ops, his criterion says that a given $P$\ has the \tp\ \ifoi
\begin{equation}
\label{tp}
\partial_x^{\beta}\partial_{\xi}^{\alpha}p_j(x^{\prime},0,0,1)
=(-1)^{j+|\alpha|}
\partial_x^{\beta}\partial_{\xi}^{\alpha}p_j(x^{\prime},0,0,-1),
\end{equation}
for all relevant $x^{\prime}$\ and for all indices $\alpha,$\ $\beta\geq 0$, and $j\leq 0$, with $p_j$\ denoting the degree-$j$\ positively 
homogeneous component (smooth for $\xi\neq 0$) in the asymptotic expansion of the symbol of $P$\ with respect to coordinates for which the boundary is given
by $x_n=0.$\ This condition is invariant under coordinate changes that preserve the set $\{x_n=0\}$, as can be proven using the standard rules of \psd\ 
calculus. Let us say that a given $p_0\in\cisx$\ satisfies the transmission condition if its positively homogeneous extension to the cotangent bundle 
(smooth except on 
the zero section) satisfies (\ref{tp}) for $j=0$. It is much easier to see (it takes only the chain rule) that this definition is also invariant under the 
appropriate coordinates changes. It is part of the content of Theorem~1 in \cite{RS}, 2.3.3.1, that, if $p\in\cisx$\ satisfies the transmission condition,
then there exists a \psd\ \op\ with the \tp\ on $X$\ whose principal symbol is $p$.

Now we prove that condition (\ref{tp}) survives the norm closure, but only for $\alpha=\beta=0$\ and $j=0$.

\begin{thm} 
\label{stil}
The kernel of $\sig:\ac\to C(\sx)$\ is equal to $\gc$. The image of $\sig$\ consists of all those functions
in $C(\sx)$\ which, over each point of the boundary, take the same value at the two covectors that vanish on the tangent space of $\bx$.
\end{thm}
\pf It follows immediately from Lemma~\ref{novogoh}\ that $\ker\sig\subseteq\gc$. On the other hand, from $\gs=\ker\sigma$, we get: 
$\gc\subseteq\overline{\ker\sigma}\subseteq\ker\sig$. This proves the first statement. 

It follows from the Stone-Weierstrass Theorem that the set of all $p\in\cisx$\ satisfying the transmission condition (which is contained in
the image of $\sigma$, by the preceding remarks) is dense in $C(\sxt)$, with $\sxt$\ defined after (\ref{cmpsq}). 
Since the image of $\sig$\ is closed, that finishes the proof. \cqd

\begin{cor}
The principal symbol homomorphism induces a \cst-algebra isomorphism between $\ac/\gc$\ and $C(\sxt)$, 
\end{cor}

The following generalization of Gohberg and Seeley's norm estimate follows immediately from Theorem~\ref{stil}\ and the fact that
every injective \cst-algebra homomorphism is an isometry.

\begin{cor}
\label{gengo} If $P$\ is a classical zero-order \psd\ \op\ with the \tp\ on $X$, then 
$
\inf_{G}||\pp+G||=||\sigma(P)||,
$
where the infimum is taken over all polyhomogeneous zero-order singular Green \ops.
\end{cor}

The next theorem is not new (see \cite{RS}, 2.3.4.4, Corollary~2), but it is perhaps appropriate to offer here this proof. 

\begin{thm}
\label{ontogreen}
The boundary principal symbol induces an isomorphism from $\gc/\kc$\ to $C(\sbx)\otimes\kc_{\rp}$, where $\kc_{\rp}$\ denotes the ideal of compact \ops\ 
on $L^2(\rp)$.
\end{thm}

\pf That $\gam$\ defines an isometry $\gc/\kc\to C(\sbx)\otimes\kc_{\rp}$\ follows from (\ref{rs}), since $\sigma(P)=0$\ if $P$\ has negative order,
and $A_{jk}$\ is compact, if $j\neq 1$\ or $k\neq 1$, for all $A\in\ac$. To prove that this isometry is surjective, it is enough to show that $\img$\ contains 
$C_0(S^{*}V,\kc_{\rp})$, for any $V$\ as in Lemma~\ref{trivial}. 

Given $f,g\in\srp$\ and $p\in C_c^\infty(S^{*}V)$, let $\tilde p(\cvb)$\ denote the local expression of the homogeneous extension of $p$\ of degree
$(-1)$.
Denoting also by $f$\ and $g$\ their extensions to $\rp$\ vanishing on the negative half-axis, let $\varphi$\ and $\psi$\ denote the Fourier transforms 
of $g$\ and $h$, respectively, where $h(t)=f(-t)$. We then define
$
\tilde g(\cvb,\xi_n,\eta_n)=\tilde p(\cvb)\omega(|\xi^\prime|)\varphi(\frac{\xi_n}{|\xi^\prime|})\psi(\frac{\eta_n}{|\xi^\prime|}),
$\ 
with $\omega$, as in Section~\ref{kg}, denoting an excision function.

It is straightforward to check (using \cite{G}, (1.2.38) and (2.3.25), for example) that $\tilde g$\ is a singular Green symbol of order and class zero on 
the euclidean space. Moreover, it is polyhomogeneous and its homogeneous principal part (smooth for $\xi^\prime\neq 0$) is given by
$
\tilde g_0(\cvb,\xi_n,\eta_n)=\tilde p(\cvb)\varphi(\frac{\xi_n}{|\xi^\prime|})\psi(\frac{\eta_n}{|\xi^\prime|}).
$ 
It then follows (from \cite{G}, (2.3.25) and (2.4.6), for example) that 
\begin{equation}
\label{derradeira}
\tilde g_0(\cvb,D_n)=|\xi^\prime|\,\tilde p(\cvb)\cdot(\kappa_{|\xi^\prime|}\circ|g\rangle\circ\langle f|\circ\kappa_{|\xi^\prime|}^{-1}), 
\end{equation}
with $\circ$\ denoting, above, composition of \ops\ on $L^2(\rp)$.

Let $G$\ denote the pullback to $X$\ of $\rho_1\tilde G\rho_1$, with $\tilde G$\ denoting the singular Green \op\ of symbol $\tilde g$\ and $\rho_1$\ as in the
proof of Lemma~\ref{trivial}. It follows from (\ref{derradeira}) that $\gamma(G)=\iota(p\otimes(|g\rangle\langle f|))$, with $\iota$\ denoting the Banach 
space isomorphism of $C_0(S^{*}V,\kc_{\rp})$\ onto itself
\[
F(\cvb)\mapsto |\xi^\prime|\cdot(\kappa_{|\xi^\prime|}\circ F(\cvb)\circ\kappa_{|\xi^\prime|}^{-1}).
\]
We are finished, because the set of all such $p\otimes(|g\rangle\langle f|)$\ is dense. \cqd

>From Theorems~\ref{stil}\ and~\ref{ontogreen}, we obtain the composition sequence $0\subset\kc\subset\gc\subset\ac$, with $\ac/\gc\simeq C(\sxt)$\
and $\gc/\kc\simeq C(\sbx)\otimes\kc_{\rp}$, isomorphisms induced by the principal symbol and by the boundary principal symbol, respectively. 

\section*{Acknowledgements}

Conversations or email with many friends and colleagues had an influence, sometimes decisive, at various stages of this work. It is a pleasure for us 
to thank them all, and to name some of them: Severino Collier Coutinho, Ruy Exel, Robert Lauter, Sergiu Moroianu, Anton Savin, Thomas Schick, Mariusz 
Wodzicki.

This work was done while the first named author was a guest of the Institute of Mathematics of the University of Potsdam. He wishes to express his gratitude to
all members of the group Partial Differential Equations and Complex Analysis for their always gentle hospitality. His visit was supported by the Funda\c c\~ao 
de Amparo \`a Pesquisa do Estado de S\~ao Paulo (Fapesp, Brazil, Processo 00/00451-2), the European research and training network ``Geometric Analysis'' (Contract HPRN-CT-1999-00118), and the Deutsche Forschungsgemeinschaft (DFG, Germany).

{\footnotesize 

\vskip0.12cm

2000 Mathematics Subject Classification: 58J40, 46L80 (58J32, 35S15, 19K56).

\vskip0.12cm
 
Instituto de  Matem\'atica e Estat\'\i stica, Universidade de S\~ao Paulo,

Caixa Postal 66281, 05315-970 S\~ao Paulo, Brazil.

melo@ime.usp.br

\vskip0.08cm

Department of Mathematics, University of Copenhagen,

Universitetsparken 5, 2100 Copenhagen, Denmark.

rnest@math.ku.dk

\vskip0.08cm

Universit\"at Potsdam, Institut f\"ur Mathematik, 

Postfach 601553, 14415 Potsdam, Germany.

schrohe@math.uni-potsdam.de

}

\end{document}